\documentclass[11pt,doublespace]{article}
\begin{document}


\newtheorem{thm}{Theorem}[section]


\newtheorem{lem}[thm]{Lemma}
\newtheorem{cor}[thm]{Corollary}
\newtheorem{ex}[thm]{Example}
\newtheorem{prop}[thm]{Proposition}
\newtheorem{remark}[thm]{Remark}
\newtheorem{coun}[thm]{Counterexample}
\newtheorem{defn}[thm]{Definition}
\newtheorem{conj}{Conjecture}
\newtheorem{problem}{Problem}
\newcommand\ack{\section*{Acknowledgement.}}

\newcommand{\etal}{{\it et al. }}

\newcommand{\bbP}{{\rm I\hspace{-0.8mm}P}}
\newcommand{\bbE}{{\rm I\hspace{-0.8mm}E}}
\newcommand{\bbF}{{\rm I\hspace{-0.8mm}F}}
\newcommand{\bbI}{{\rm I\hspace{-0.8mm}I}}
\newcommand{\bbR}{{\rm I\hspace{-0.8mm}R}}
\newcommand{\bbRp}{{\rm I\hspace{-0.8mm}R}_+}
\newcommand{\bbN}{{\rm I\hspace{-0.8mm}N}}
\newcommand{\bbC}{{\rm C\hspace{-2.2mm}|\hspace{1.2mm}}}
\newcommand{\bbD}{{\rm I\hspace{-0.8mm}D}}
\newcommand{\bbQ}{\bf Q}
\newcommand{\bbZ}{{\rm \rlap Z\kern 2.2pt Z}}
\newcommand{\bbK}{{\rm I\hspace{-0.8mm}K}}

\newcommand{\matP}{{\bbP}}
\newcommand{\mattildeP}{\tilde{\bbP}}
\newcommand{\matPN}[1]{{\bbP}_{#1}}
\newcommand{\matPP}[1]{{\bbP}_{#1}^0}
\newcommand{\matE}{{\bbE}}
\newcommand{\mattildeE}{\tilde{\bbE}}
\newcommand{\matEP}[1]{{\bbE}_{#1}^0}
\newcommand{\matF}{{\bbF}}
\newcommand{\matR}{{\bbR}}
\newcommand{\matRp}{{\bbRp}}
\newcommand{\matN}{{\bbN}}
\newcommand{\matZ}{{\bbZ}}
\newcommand{\matI}{{\bbI}}
\newcommand{\matK}{{\bbK}}
\newcommand{\matQ}{{\bbQ}}
\newcommand{\matC}{{\bbC}}
\newcommand{\matD}{{\bbD}}

\newcommand{\calL}{{\cal L}}
\newcommand{\calM}{{\cal M}}
\newcommand{\calN}{{\cal N}}
\newcommand{\calF}{{\cal F}}
\newcommand{\calG}{{\cal G}}
\newcommand{\calD}{{\cal D}}
\newcommand{\calB}{{\cal B}}
\newcommand{\calH}{{\cal H}}
\newcommand{\calI}{{\cal I}}
\newcommand{\calP}{{\cal P}}
\newcommand{\calQ}{{\cal Q}}
\newcommand{\calS}{{\cal S}}
\newcommand{\calT}{{\cal T}}
\newcommand{\calC}{{\cal C}}
\newcommand{\calK}{{\cal K}}
\newcommand{\calX}{{\cal X}}
\newcommand{\cals}{{\cal S}}
\newcommand{\calE}{{\cal E}}

\newcommand{\koniecmat}{\,}

\newcommand{\eqd}{\ =_{\rm d}\ }
\newcommand{\toto}{\leftrightarrow}
\newcommand{\eqdistr}{\stackrel{\rm d}{=}}
\newcommand{\as}{\stackrel{\rm a.s.}{=}}
\newcommand{\assubs}{{\rm a.s.}}
\newcommand{\convdistr}{\stackrel{\rm d}{\rightarrow}}
\newcommand{\convweak}{{\Rightarrow}}
\newcommand{\convas}{\stackrel{\rm a.s.}{\rightarrow}}
\newcommand{\convfidi}{\stackrel{\rm fidi}{\rightarrow}}
\newcommand{\convprob}{\stackrel{p}{\rightarrow}}
\newcommand{\deff}{\stackrel{\rm def}{=}}
\newcommand{\bis}{{'}{'}}
\newcommand{\Cov}{{\rm Cov}}
\newcommand{\Var}{{\rm Var}}
\newcommand{\Exp}{{\rm E}}

\newcommand{\nd}{n^{\delta}}
\newcommand{\koniec}{\newline\vspace{3mm}\hfill $\odot$}


\title{Reduction principles for quantile and Bahadur-Kiefer processes of long-range dependent linear sequences  \protect}
\author{Mikl\'{o}s Cs\"{o}rg\H{o}\thanks{School of Mathematics and Statistics, Carleton University, 1125 Colonel By Drive, Ottawa, Ontario, K1S 5B6 Canada, email: mcsorgo@math.carleton.ca
} \and Rafa{\l} Kulik\thanks{School of Mathematics and Statistics,
University of Sydney, NSW 2006, Australia, email:
rkuli@maths.usyd.edu.au and Mathematical Institute, Wroc{\l}aw
University, Pl. Grunwaldzki 2/4, 50-384 Wroc{\l}aw, Poland.
\newline Research supported in part by NSERC Canada Discovery Grants of
Mikl\'{o}s Cs\"{o}rg\H{o}, Donald Dawson and Barbara Szyszkowicz at
Carleton University}}

\maketitle

\begin{center}{\small
Title of the document: Reduction principles LRD\\

The final version (before editorial proofs) accepted to Probability
Theory and Related Fields on 15 October 2007.}
\end{center}

\begin{abstract}
In this paper we consider quantile and Bahadur-Kiefer processes for
long range dependent linear sequences. These processes, unlike in
previous studies, are considered on the whole interval $(0,1)$. As
it is well-known, quantile processes can have very erratic behavior
on the tails. We overcome this problem by considering these
processes with appropriate weight functions. In this way we conclude
strong approximations that yield some remarkable phenomena that are
not shared with i.i.d. sequences, including weak convergence of the
Bahadur-Kiefer processes, a different pointwise behavior of the
general and uniform Bahadur-Kiefer processes, and a somewhat
"strange" behavior of the general quantile process.
\end{abstract}

\noindent {\bf Short title:} Quantiles and LRD

\section{Introduction}
Let $\{\epsilon_i,i\ge 1\}$ be a centered sequence of i.i.d. random
variables. Consider the class of stationary linear processes
\begin{equation}\label{model}
X_i=\sum_{k=0}^{\infty}c_k\epsilon_{i-k},\ \ \ i\ge 1 .
\end{equation}
We assume that the sequence $c_k$, $k\ge 0$, is regularly varying
with index $-\beta$, $\beta\in (1/2,1)$ (written as $c_k\in
RV_{-\beta}$). This means that $c_k\sim k^{-\beta}L_0(k)$ as
$k\to\infty$, where $L_0$ is slowly varying at infinity. We shall
refer to all such models as long range dependent (LRD) linear
processes. In particular, if the variance exists, then the
covariances $\rho_k:=\Exp X_0X_k$ decay at the hyperbolic rate,
$\rho_k=L(k)k^{-(2\beta-1)}=:L(k)k^{-D}$, where
$\lim_{k\to\infty}L(k)/L_0^2(k)=B(2\beta-1,1-\beta)$ and
$B(\cdot,\cdot)$ is the beta-function. Consequently, the covariances
are not summable (cf. \cite{GiraitisSurgailis2002}).

Assume that $X_1$ has a continuous distribution function $F$. For
$y\in (0,1)$ define $Q(y)=\inf\{x:F(x)\ge y\}=\inf\{x:F(x)= y\}$,
the corresponding (continuous) quantile function. Given the ordered
sample $X_{1:n}\le\cdots\le X_{n:n}$ of $X_1,\ldots,X_n$, let
$F_n(x)=n^{-1}\sum_{i=1}^n1_{\{X_i\le x\}}$ be the empirical
distribution function and $Q_n(\cdot)$ be the corresponding
left-continuous sample quantile function. Define $U_i=F(X_i)$ and
$E_n(x)=n^{-1}\sum_{i=1}^n1_{\{U_i\le x\}}$, the associated uniform
empirical distribution. Denote by $U_n(\cdot)$ the corresponding
uniform sample quantile function.

Our purpose in this paper is to study the asymptotic behavior of
sample quantiles for long range dependent sequences. This will be
done in the spirit of the Bahadur-Kiefer approach (cf.
\cite{Bahadur1966}, \cite{Kiefer1967}, \cite{Kiefer1970}).

Assume that $E\epsilon_1^2<\infty$. Let $r$ be an integer and define
$$Y_{n,r}=\sum_{i=1}^n\sum_{1\le j_1<\cdots\le
j_r}\prod_{s=1}^rc_{j_s}\epsilon_{i-j_s},\qquad n\ge 1,$$ so that
$Y_{n,0}=n$, and $Y_{n,1}=\sum_{i=1}^nX_i$. If $p<(2\beta-1)^{-1}$,
then
\begin{equation}\label{eq-varance-behaviour}
\sigma_{n,p}^2:=\Var (Y_{n,p})\sim n^{2-p(2\beta-1)}L_0^{2p}(n).
\end{equation}
Define now the general empirical, the uniform empirical, the general
quantile and the uniform quantile processes respectively as follows:
\begin{equation}\label{eq.general.empirical}
\beta_n(x)=\sigma_{n,1}^{-1}n(F_{n}(x)-F(x)),\qquad x\in {\matR},
\end{equation}
\begin{equation}\label{eq.uniform.empirical}
\alpha_n(y)=\sigma_{n,1}^{-1}n(E_{n}(y)-y),\qquad y\in (0,1),
\end{equation}
\begin{equation}\label{eq.general.quantile}
q_n(y)=\sigma_{n,1}^{-1}n(Q(y)-Q_{n}(y)),\qquad y\in (0,1),
\end{equation}
\begin{equation}\label{eq.uniform.quantile}
u_n(y)=\sigma_{n,1}^{-1}n(y-U_{n}(y)),\qquad y\in (0,1).
\end{equation}
Assume for a while that $X_i$, $i\ge 1$ are i.i.d. We shall refer to
this as to the {\it i.i.d. model}. Denote by $\alpha_n^{\mbox{\rm
iid}}$, $q_n^{\mbox{\rm iid}}$, $u_n^{\mbox{\rm iid}}$ the uniform
empirical, general quantile, uniform quantile processes based on
i.i.d. samples with the constants $\sigma_{n,1}^{-1}n$ in
(\ref{eq.uniform.empirical}), (\ref{eq.general.quantile}),
(\ref{eq.uniform.quantile}) replaced with $\sqrt{n}$. Fix $y\in
(0,1)$. Let $I_y$ be a neighborhood of $Q(y)$ and assume that $F$ is
twice differentiable with respect to Lebesgue measure with
respective first and second derivatives $f$ and $f'$. Assuming that
$\inf_{x\in I_y}f(x)>0$ and $\sup_{x\in I_y}|f^{'}(x)|<\infty$,
Bahadur in \cite{Bahadur1966} obtained the following {\it Bahadur
representation} of quantiles
\begin{equation}\label{Bahadur}
\alpha_n^{\mbox{\rm iid}}(y)-f(Q(y))q_n^{\mbox{\rm
iid}}(y)=:R_n^{\mbox{\rm iid}}(y),
\end{equation}
with
\begin{equation}\label{Bahadur-rate}
R_n^{\mbox{\rm iid}}(y)=O_{\assubs}(n^{-1/4}(\log n)^{1/2}(\log\log
n)^{1/4}),\qquad n\to\infty ,
\end{equation}
The process $\{R_n^{\mbox{\rm iid}}(y),y\in (0,1)\}$ is called the
Bahadur-Kiefer process. Later, Kiefer proved in \cite{Kiefer1967}
that (\ref{Bahadur-rate}) can be strengthened to
\begin{equation}\label{Bahadur-rate-improvement}
R_n^{\mbox{\rm iid}}(y)=O_{\assubs}(n^{-1/4}(\log\log n)^{3/4}),
\end{equation}
which is the optimal rate. Continuing his study, in
\cite{Kiefer1970} Kiefer obtained the uniform version of
(\ref{Bahadur}), referred to later on as the {\it Bahadur-Kiefer
representation}:
\begin{equation}\label{Kiefer}
\sup_{y\in [0,1]}\left|\alpha_n^{\mbox{\rm
iid}}(y)-f(Q(y))q_n^{\mbox{\rm iid}}(y)\right|=:R_n^{\mbox{\rm iid}}
\end{equation}
where
\begin{equation}\label{Kiefer-rate}
R_n^{\mbox{\rm iid}}=O_{\assubs}(n^{-1/4}(\log n)^{1/2}(\log\log
n)^{1/4}),\qquad n\to\infty .
\end{equation}
Once again, the above rate is optimal. Kiefer obtained his result
assuming
\begin{itemize}
\item[{\rm (K1)}] $f$ has finite support and $\sup_{x\in {\matR}}|f^{'}(x)|<\infty$,
\item[{\rm (K2)}] $\inf_{x\in {\matR}}f(x)>0$.
\end{itemize}
We shall refer to (K1)-(K2) as to the {\it Kiefer conditions}.

Further on, Cs\"{o}rg\H{o} and R\'{e}v\'{e}sz
\cite{CsorgoRevesz1978} obtained Kiefer's result (\ref{Kiefer})
under the following, weaker conditions, which shall be referred to
later on as the {\it Cs\"{o}rg\H{o}-R\'{e}v\'{e}sz conditions} (cf.
also \cite[Theorem 3.2.1]{CsorgoLN}):
\begin{itemize}
\item[\rm{(CsR1)}] $f$ exists on $(a,b)$, where
$a=\sup\{x:F(x)=0\}$, $b=\inf\{x:F(x)=1\}$, $-\infty\le a<b\le
\infty$,
\item[{\rm (CsR2)}] $\inf_{x\in (a,b)}f(x)>0$,
\item[{\rm (CsR3)}] $\sup_{x\in (a,b)}F(x)(1-F(x))\frac{|f'(x)|}{f^2(x)}=\sup_{y\in
(0,1)}y(1-y)\left|\frac{f'(Q(y))}{f^2(Q(y))}\right|\le \gamma$ with
some $\gamma>0$,
\item[{\rm (CsR4)}] (i) $0<A:=\lim_{y\downarrow 0}f(Q(y))<\infty$, $0<B:=\lim_{y\uparrow
1}f(Q(y))<\infty$, or\\
(ii) if $A=0$ (respectively $B=0$) then $f$ is nondecreasing
(respectively nonincreasing) on an interval to the right of $Q(0+)$
(respectively to the left of $Q(1-)$).
\end{itemize}
In particular, they showed that, under (CsR1), (CsR2), (CsR3), as
$n\to\infty$,
\begin{equation}\label{CR-1}
\sup_{n^{-1}\log\log n\le y\le 1-n^{-1}\log\log
n}|f(Q(y))q_n^{\mbox{\rm iid}}(y)-u_n^{\mbox{\rm
iid}}(y)|=O_{\assubs}(n^{-1/2}\log\log n).
\end{equation}
Additionally, if (CsR4) holds, then, as $n\to\infty$,
\begin{equation}\label{CR-2}
\sup_{y\in [0,1]}|f(Q(y))q_n^{\mbox{\rm iid}}(y)-u_n^{\mbox{\rm
iid}}(y)|=O_{\assubs}(n^{-1/2}\ell (n)).
\end{equation}
Here, and in the sequel, $\ell(n)$ is a slowly varying function at
infinity, but can be different at each place it appears (e.g. when
Cs\"{o}rg\H{o}-R\'{e}v\'{e}sz conditions hold, then
$\ell(n)=\log\log n$). This, via the special case of
(\ref{Kiefer-rate})
$$ \sup_{y\in
[0,1]}|u_n^{\mbox{\rm iid}}(y)-\alpha_n^{\mbox{\rm
iid}}(y)|=O_{\assubs}(n^{-1/4}(\log n)^{1/2}(\log\log n)^{1/4}),
$$
yields the Bahadur-Kiefer representation (\ref{Kiefer-rate}) under
less restrictive conditions compared to Kiefer's assumptions. In
particular, Cs\"{o}rg\H{o}-R\'{e}v\'{e}sz conditions are fulfilled
if $F$ is exponential or normal. Also, if (CsR4)(i) obtains, then
$\ell (n)$ in (\ref{CR-2}) is $\log\log n$. We refer to
\cite{CsorgoLN}, \cite{CsorgoReveszbook} and
\cite{CsorgoSzyszkowicz1998} for more discussion of these
conditions. We note in passing that taking $\sup$ over
$[1/(n+1),n/(n+1)]$ instead of the whole unit interval, the
statement (\ref{CR-2}) holds true assuming only the conditions
(CsR1)-(CsR3) (cf. \cite[Theorem 3.1]{CsorgoCsorgoHorvathRevesz}, or
\cite[Theorem 6.3.1]{CsorgoHorvath-book1993}).

As to LRD linear processes with partial sums $Y_{n,r}$ above, the
first result on sample quantiles can be found in Ho and Hsing
\cite{HoHsing}, where it is shown under Kiefer-type conditions that,
as $n\to\infty$, one has for all $\beta\in (\frac{1}{2},1)$
\begin{equation}\label{quantiles-HoHsing}
\sup_{y\in
(y_0,y_1)}|Q(y)-Q_n(y)-n^{-1}Y_{n,1}|=o_{\assubs}(n^{-(1+\lambda)}\sigma_{n,1}),
\end{equation}
where $0<y_0<y_1<1$ are fixed and
$0<\lambda<(\beta-\frac{1}{2})\wedge (1-\beta)$. This means that the
sample quantiles $Q_n(y)$, $y\in (y_0,y_1)$ can be approximated by
the sample mean $n^{-1}Y_{n,1}=n^{-1}\sum_{i=1}^nX_i$ independently of $y$. This quantile
process approximation is a consequence of their landmark result for
empirical processes; see also \cite{Surgailis-Koul-2002},
\cite{Surgailis-2002} and \cite{Surgailis-2004} for related studies.
The best available result along these lines is due to Wu
\cite{Wu2003}. To state a particular version of his result, let
$F_{\epsilon}$ be the distribution function of the centered i.i.d.
sequence $\{\epsilon_i,i\ge 1\}$. Assume that for a given integer
$p$, the derivatives $F^{(1)}_{\epsilon}, \ldots,
F^{(p+3)}_{\epsilon}$ of $F_{\epsilon}$ are bounded and integrable.
Note that these properties are inherited by the distribution $F$ as
well (cf. \cite{HoHsing} or \cite{Wu2003}).
\begin{thm}\label{thm-HoHsing}
Let $p$ be a positive integer. Then, as $n\to\infty$,
$$
\Exp \sup_{x\in {\matR}}\left|\sum_{i=1}^n(1_{\{X_i\le
x\}}-F(x))+\sum_{r=1}^p(-1)^{r-1}F^{(r)}(x)Y_{n,r}\right|^2=O(\Xi_n+n(\log
n)^2),
$$
where
$$
\Xi_n=\left\{\begin{array}{ll} O(n), & (p+1)(2\beta-1)>1\\
O(n^{2-(p+1)(2\beta-1)}L_0^{2(p+1)}(n)), & (p+1)(2\beta-1)<1
\end{array}\right. .
$$
\end{thm}

Using this result, under Kiefer conditions as $n\to\infty$, Wu
\cite{Wu2005} obtained
\begin{equation}\label{BahadurKiefer-Wu}
\sup_{y\in
(y_0,y_1)}|\alpha_n(y)-f(Q(y))q_n(y)-\sigma_{n,1}^{-1}n^{-1}Y_{n,1}^2f^{'}(Q(y))/2|=O_{\assubs}(j_n\ell
(n)),
\end{equation}
where $j_n=n^{-(\frac{3}{4}-\frac{\beta}{2})}$ if
$\beta>\frac{7}{10}$ and $j_n=n^{-(2\beta-1)}$ if
$\beta\le\frac{7}{10}$. As argued in \cite[Section 7.1]{Wu2005} this
bound is sharp up to a multiplicative slowly varying function $\ell
(n)$. From (\ref{BahadurKiefer-Wu}) and the central limit theorem for the partial sums $\sum_{i=1}^nX_i$ we may also deduce under Kiefer
conditions and $\beta\in (\frac{1}{2},\frac{5}{6})$, that for the
Bahadur-Kiefer process
\begin{equation}\label{BK-definition}
R_n(y)=\alpha_n(y)-f(Q(y))q_n(y)
\end{equation}
we have weak convergence $\sigma_{n,1}^{-1}nR_n(y)\convweak
f'(Q(y))Z^2/2$ in $D([y_0,y_1])$,Cs\"{o}rg\H{o}-R\'{e}v\'{e}sz
conditions equipped with the sup-norm topology, where $Z$ is a
standard normal random variable. In particular, if $\epsilon_i, i\ge
1$ are i.i.d. standard normal random variables, then, as
$n\to\infty$,
\begin{equation}\label{eq-BahadurKiefer-general-normal}
\sigma_{n,1}^{-1}nR_n(y)\convweak \phi'(\Phi^{-1}(y))Z^2/2 \qquad
\mbox{\rm in }D([y_0,y_1]),
\end{equation}
where $\phi$ and $\Phi$ are the standard normal density and
distribution functions, respectively.

This behavior is completely different compared to the i.i.d. case,
for it is well known that the Bahadur-Kiefer process cannot converge
weakly in the space of cadlag functions (cf., e.g., \cite[Remark
2.1]{CSW2006}).

However, this weak convergence phenomenon was first observed
explicitly by Cs\"{o}rg\H{o}, Szyszkowicz and Wang \cite{CSW2006}
for long range dependent Gaussian sequences. For the sake of
comparison with (\ref{eq-BahadurKiefer-general-normal}), assume that
$\epsilon_i,i\ge 1$ are standard normal random variables and that
$\sum_{k=1}^{\infty}c_k^2=1$. Then the $X_i$ defined by
(\ref{model}) are standard normal. Define $Y_n=G(X_n)$, with some
real-valued measurable function $G$. Let $J_l(y)=\Exp
\left[\left(1_{\{F(G(X))\le y\}}-y\right)H_l(X)\right]$, where $H_l$
is the $l$th Hermite polynomial. In particular, taking
$G=F^{-1}\Phi$ we have that $Y_n$ have the marginal distribution
$F$. The Hermite rank is 1 and $J_1(y)=-\phi(\Phi^{-1}(y))$, and we
may take $Y_n=X_n$. Note that for the Hermite rank 1, via $L(n)\sim
B(2\beta-1,1-\beta)L^2_0(n)$, their scaling factor $d_n^2=n^{2-\tau
D}L^{\tau}(n)$ (cf. (1.5) of \cite{CSW2006}) agrees (up to a
constant) with $\sigma_{n,1}^2$ of (\ref{eq-varance-behaviour}).
Note also that $J_1(y)J_1'(y)=\phi'(\Phi^{-1}(y))$. Thus, for the
uniform Bahadur-Kiefer process
\begin{equation}\label{eq-definition-BahadurKiefer-uniform}
\tilde{R}_n(y)=\alpha_n(y)-u_n(y)
\end{equation}
we may conclude from \cite[Theorem 2.3]{CSW2006} that (see also
Remark \ref{rem-gaussian} in the present paper), as $n\to\infty$,
\begin{equation}\label{eq-uniform-BahadurKiefer-gaussian}
\sigma_{n,1}^{-1}n\tilde R_n(y)\convweak \phi'(\Phi^{-1}(y))Z^2
\qquad \mbox{\rm in }D([y_0,y_1]).
\end{equation}
Comparing (\ref{eq-BahadurKiefer-general-normal}) with
(\ref{eq-uniform-BahadurKiefer-gaussian}), we see that the weak
limits in $D[y_0,y_1]$ of the uniform and the general Bahadur-Kiefer
processes are different.

We note that Cs\"{o}rg\H{o} \etal \cite{CSW2006} have also established
the rate for the deviation of $\tilde R_n(y)$ from $R_n(y)$ under
the Cs\"{o}rg\H{o}-R\'{e}v\'{e}sz conditions. This rate, in the case
of the Hermite rank 1, coincides with the scaling factor for the
weak convergence of the Bahadur-Kiefer processes in
(\ref{eq-BahadurKiefer-general-normal}) and
(\ref{eq-uniform-BahadurKiefer-gaussian}). Since the uniform and the
general Bahadur-Kiefer processes have different limits, the
rate obtained for their nearness in \cite{CSW2006} cannot be improved. \\

In this paper we deal with several problems. First, unlike in
\cite{HoHsing} or \cite{Wu2005}, we consider quantile and
Bahadur-Kiefer processes on the whole interval $(0,1)$ under very
general conditions on the distribution function $F$. As it is
well-known, quantile processes can have very erratic behavior on the
tails. Moreover, it should be pointed out that in the LRD case, even
when we deal with the associated uniform version of quantile and
Bahadur-Kiefer processes, we also have to deal with the general
quantile function of $X_1$. We solve this problem by considering
these processes with appropriate weight functions. With this help,
we can conclude various strong approximations, as well as some
remarkable phenomena not shared with i.i.d. sequences, including
weak convergence of the Bahadur-Kiefer processes, or different
pointwise behavior of the general and uniform Bahadur-Kiefer
processes. Further on, we deal with the general quantile process
$q_n(y)$. Via its weak convergence, we obtain confidence intervals
for the quantile function $Q$. Moreover, if one considers the
subordinated Gaussian sequence $Y_n=G(X_n)$, then the behavior of
the quantile process does not only depend on the marginals of
$Y_n$'s and the dependence structure (i.e. the parameter $\beta$),
but also on a "hidden" LRD sequence $\{X_i,i\ge 1\}$. This property
cannot occur in a weakly dependent case.

Although, especially by dealing with weight functions, the paper is
fairly technical, however, the choice of 'good' weight functions
allow us to obtain reasonable simultaneous confidence intervals for
the quantile function (see Section \ref{sec-weak-behaviour}).

Our results are presented in Section \ref{sec-results}. That section
is concluded with a number of remarks (see Section
\ref{sec-remarks}), including a discussion of the recent paper \cite{CSW2006}. The proofs are given in Section \ref{sec-proofs}\\

In what follows $C$ will denote a generic constant which may be
different at each of its appearances. Also, for any sequences $a_n$
and $b_n$, we write $a_n\sim b_n$ if $\lim_{n\to\infty}a_n/b_n=1$.
Further, recall that $\ell(n)$ is a slowly varying function,
possibly different at each place it appears. Moreover, $f^{(k)}$ denotes the $k$th order derivative of $f$.\\
\section{Statement of results and discussion}\label{sec-results}
For discussing our results, we introduce some notation.

\noindent Let $p$ be a positive integer and put
\begin{eqnarray*}
S_{n,p}(x) &=&\sum_{i=1}^n(1_{\{X_i\le
x\}}-F(x))+\sum_{r=1}^p(-1)^{r-1}F^{(r)}(x)Y_{n,r}\\
&=:&\sum_{i=1}^n(1_{\{X_i\le x\}}-F(x))+V_{n,p}(x),
\end{eqnarray*}
so that $ S_{n,1}(x)=nF_n(x)+f(x)\sum_{i=1}^nX_i , $ and $
S_{n,0}(x)=nF_n(x) . $ Setting $U_i=F(X_i)$ and $x=Q(y)$ in the
definition of $S_{n,p}(\cdot)$, we arrive at its uniform version,
\begin{eqnarray}\label{eq.uniform.version}
\tilde S_{n,p}(y)&=&\sum_{i=1}^n(1_{\{U_i\le
y\}}-y)+\sum_{r=1}^p(-1)^{r-1}F^{(r)}(Q(y))Y_{n,r}\nonumber\\
&=:&\sum_{i=1}^n(1_{\{U_i\le y\}}-y)-\tilde{V}_{n,p}(y) .
\end{eqnarray}
Recall that
$$
R_n(y)=\alpha_n(y)-f(Q(y))q_n(y),\qquad y\in (0,1),
$$
is the Bahadur-Kiefer process and
$$
\tilde{R}_n(y)=\alpha_n(y)-u_n(y),\qquad y\in (0,1),
$$
is the uniform Bahadur-Kiefer process.\\

We shall consider the following assumptions on the distribution
function $F$.
\begin{itemize}
\item[{\bf (A(p))}] The functions $(f^{(r-1)}\circ Q)^{(1)}(y)$,
$r=1,\ldots,p$, are uniformly bounded. The integer $p$ will be
chosen appropriately in the sequel.
\item[{\bf (B)}] The function $(f\circ Q)^{(2)}(y)$ is
uniformly bounded.
\item[{\bf (C(p))}] For $r=0,\ldots,p-1$,
$$
\sup_{y\in (0,1)}\frac{f^{(r+1)}(Q(y))}{f(Q(y))}(y(1-y))^{1/2}=O(1).
$$
\end{itemize}
\subsection{Strong approximations}
Let
$$
a_n=\sigma_{n,1}n^{-1}\log\log
n=n^{-(\beta-\frac{1}{2})}L_0(n)\log\log n ,
$$
$$
b_n=\sigma_{n,1}^2n^{-1}a_n(\log\log
n)^{1/2}=n^{-(3\beta-\frac{5}{2})}L_0^3(n)(\log\log n)^{3/2},
$$
$$
c_n=\sigma_{n,1}^{-1}b_n(\log
n)^{1/2}=n^{-(2\beta-1)}L_0^2(n)(\log\log n)^{3/2}(\log n)^{1/2},
$$
$$
d_{n,p}=\left\{\begin{array}{ll} n^{-(1-\beta)}L_0^{-1}(n)(\log n)^{5/2}(\log\log n)^{3/4}, & (p+1)(2\beta-1)>1\\
n^{-p(\beta-\frac{1}{2})}L_0^{p}(n)(\log n)^{1/2}(\log\log n)^{3/4},
& (p+1)(2\beta-1)<1
\end{array}\right. ,
$$
$$
b_{n,p}=\sigma_{n,1}^2n^{-1}d_{n,p}(\log\log n)^{1/2},
$$
and
$$\delta_n=n^{-(2\beta-1)}L_0^2(n)(\log\log n).
$$
\subsubsection{Reduction principles for the uniform quantile process}
First, we deal with reduction principles for quantiles. Ho and
Hsing, \cite[p. 1003]{HoHsing} asked, whether there was an expansion
for the quantile process which mirrors that in their Theorem 2.1 for
the empirical process. We have the following result.
\begin{thm}\label{cor-reduction-quantiles}
Assume {\rm (B)}, and either {\rm (A(1))} or {\rm (A(2))} according
to $\beta\ge 3/4$ or $\beta<3/4$. Then, under the conditions of {\rm
Theorem {\rm\ref{thm-HoHsing}}}, as $n\to\infty$, we have,
\begin{equation}\label{eq-uniform-quantile-approximation}
\sup_{y\in
(0,1)}\left|u_n(y)+\sigma_{n,1}^{-1}f(Q(y))\sum_{i=1}^nX_i\right|=\left\{
\begin{array}{ll}
O_{\assubs}(d_{n,1}), &\mbox{\rm if}\;\; \beta\ge 3/4,\\
O_{a.s.}(a_n), &\mbox{\rm if}\;\; \beta<3/4  .
\end{array}
\right.
\end{equation}
If $\beta<\frac{3}{4}$, the bound is is optimal.
\end{thm}

To remove assumptions (A) and (B) we shall consider a (possibly)
{\it weighted} approximation of the uniform quantiles. Define
$\psi_1(y)$ in the following way. If $\beta<\frac{3}{4}$, then
$\psi_1(y)=1$ if (C(2)) holds, and
$\psi_1(y)=(y(1-y))^{\gamma-\frac{1}{2}+\mu}$, $\mu>0$ otherwise. If
$\beta\ge \frac{3}{4}$, $\psi_1(y)=(y(1-y))^{\gamma+\mu}$.
\begin{thm}\label{th-reduction-quantiles-general}
Let $p=2$. Then, under the conditions of {\rm Theorem}
{\rm\ref{thm-HoHsing}}, as $n\to\infty$, we have,
\begin{equation}\label{quantile-approx-weighted}
\sup_{y\in
(0,1)}\psi_1(y)\left|u_n(y)+\sigma_{n,1}^{-1}f(Q(y))\sum_{i=1}^nX_i\right|=\left\{
\begin{array}{ll}
O_{\assubs}(d_{n,1}), &\mbox{\rm if}\;\; \beta\ge 3/4,\\
O_{a.s.}(a_n), &\mbox{\rm if}\;\; \beta<3/4  .
\end{array}
\right.
\end{equation}
If $\beta<\frac{3}{4}$, the bound is optimal.
\end{thm}

\indent From Theorems \ref{cor-reduction-quantiles} or
\ref{th-reduction-quantiles-general} and Lemma \ref{lem.3} below we
have the following reduction principle for quantiles, which mirrors
that for the empirical process. In order to state the result,
redefine $\psi_1$ to be $1$ if, for a given $p$, (A(p)) holds, and
be it as before otherwise. The result is stated for
$\beta<\frac{3}{4}$, only in order to avoid the additional term
coming from $d_{n,1}$.
\begin{cor}\label{red-princ-quantiles}
$\beta<\frac{3}{4}$. Let $p\ge 1$ be an arbitrary integer such that
$p<(2\beta-1)^{-1}$. Assume that either {\rm(A(p))} and {\rm(B)}, or
{\rm(C(p))} hold. Under the conditions of {\rm Theorem}
{\rm\ref{thm-HoHsing}}, as $n\to\infty$,
\begin{eqnarray*}
\lefteqn{\sup_{y\in
(0,1)}\sigma_{n,p}^{-1}\psi_1(y)|y-U_n(y)+n^{-1}\tilde
V_{n,p}(y)|}\\
&=&O_{a.s.}(n^{-(2\beta-p\frac{2\beta-1}{2})}L_0^{2-p}(n)\log\log
n(\log n)^{1/2}).
\end{eqnarray*}
\end{cor}
\subsubsection{Approximations of the uniform Bahadur-Kiefer process}
Similarly to the uniform quantile process, in Theorem
\ref{th-uniform} we obtain strong approximation of the uniform
Bahadur-Kiefer process on the whole interval $(0,1)$ on assuming (A)
and (B).
\begin{thm}\label{th-uniform}
Assume {\rm(B)}, and either {\rm(A(2))} or {\rm(A(3))} according to
$\beta\ge 2/3$ or $\beta<2/3$. Under the conditions of {\rm Theorem}
{\rm\ref{thm-HoHsing}}, as $n\to\infty$,
\begin{equation}\label{main.result}
\sup_{y\in (0,1)}\left|\tilde R_n(y)-n^{-1}\sigma_{n,1}^{-1}
f^{(1)}(Q(y))\left(\sum_{i=1}^nX_i\right)^2\right|=\left\{
\begin{array}{ll}
O_{\assubs}(d_{n,2}), &\mbox{\rm if}\;\; \beta\ge 2/3,\\
O_{a.s.}(c_n), &\mbox{\rm if}\;\; \beta<2/3  .
\end{array}
\right.
\end{equation}
\end{thm}
To remove assumptions (A) and (B), we shall consider a
\emph{weighted} approximation of the uniform quantile and
Bahadur-Kiefer processes. Define for arbitrary $\mu>0$,
$$
\psi_2(y)=\left\{\begin{array}{lll} (y(1-y))^{1+\mu}, &&\mbox{\rm
if}\;\;
\beta<\frac{3}{4}\;\;\mbox{\rm  and } (C(3));\\
(y(1-y))^{1+\mu}, &&\mbox{\rm if}\;\; \beta<\frac{3}{4},
\gamma<\frac{3}{2}\;\;\mbox{\rm  and not}\;\; (C(3));\\
(y(1-y))^{\gamma-\frac{1}{2}+\mu} &&\mbox{\rm if}\;\;
\beta<\frac{3}{4}\;\;\mbox{\rm and } \gamma\ge\frac{3}{2};\\
(y(1-y))^{\gamma+\mu}, &&\mbox{\rm if}\;\; \beta\ge\frac{3}{4}.
\end{array}\right.
$$
\begin{thm}\label{th-reduction-quantiles-weighted}
Under the conditions of {\rm Theorem} {\rm\ref{thm-HoHsing}}, as
$n\to\infty$,
\begin{eqnarray}\label{main.result.1}
\sup_{y\in (0,1)}\psi_2(y)\left|\tilde
R_n(y)-n^{-1}\sigma_{n,1}^{-1}f^{(1)}(Q(y))\left(\sum_{i=1}^nX_i\right)^2\right|=\left\{
\begin{array}{ll}
O_{\assubs}(d_{n,2}), &\mbox{\rm if}\;\; \beta\ge 2/3,\\
O_{a.s.}(c_n), &\mbox{\rm if}\;\; \beta<2/3  .
\end{array}
\right.\nonumber.
\end{eqnarray}
\end{thm}

\indent From Theorem \ref{th-uniform} or
\ref{th-reduction-quantiles-weighted} and Lemma \ref{lem.3} below,
we obtain the reduction principle for the distance between the
uniform empirical and the uniform quantile processes, similar to
that of Corollary \ref{red-princ-quantiles}.
Further, an immediate corollary to Theorem \ref{th-uniform}, via the
LIL for partial sums $\sum_{i=1}^nX_i$ (see (\ref{step.5}) below),
is the following result.
\begin{cor}\label{cor.lil.uniform.BK}
Under the conditions of {\rm Theorem} {\rm\ref{th-uniform}}, if
$\beta<\frac{3}{4}$,
\begin{equation}\label{exact.rate}
\limsup_{n\to\infty} \sigma_{n,1}^{-1}n(\log\log n)^{-1}\sup_{y\in
(0,1)}|\tilde{R}_n(y)|\as c(\beta,1)\sup_{y\in
(0,1)}|f^{(1)}(Q(y))|,
\end{equation}
where
$c^2(\beta,p)=\left(\int_0^{\infty}x^{-\beta}(1+x)^{-\beta}dx\right)(1-\beta)^{-1}(3-2\beta)^{-1}$.
\end{cor}
\begin{cor}\label{cor-weak-convergence}
Under the conditions of {\rm Theorem} {\rm\ref{th-uniform}}, if
$\beta<\frac{3}{4}$,
$$
\sigma_{n,1}^{-1}n\tilde R_n(y)\convweak f^{(1)}(Q(y))Z^2 .
$$
\end{cor} The corresponding results can also be stated in the setting of Theorem
\ref{th-reduction-quantiles-weighted}.
\subsubsection{Approximation of the general Bahadur-Kiefer process}\label{sec-general-BK}
As for the general Bahadur-Kiefer process, a typical approach in the
i.i.d. case is to approximate the normalized quantiles
$f(Q(y))q_n(y)$ via the uniform quantiles and then use this to
generalize all results valid in the uniform case to the general one,
as described in the Introduction (cf. (\ref{CR-1}), (\ref{CR-2})).
This approach was also followed in \cite[Section 4]{CSW2006} as
well. However, this cannot work in the LRD case, for then the
uniform and general Bahadur-Kiefer processes have different limits
(cf. (\ref{eq-BahadurKiefer-general-normal}),
(\ref{eq-uniform-BahadurKiefer-gaussian})). Moreover, assumptions
(A) and (B) do not help in this case.

With arbitrary $\mu>0$, define
$$
\psi_3(y)=\left\{\begin{array}{lll} (y(1-y))^{1+\mu}, &&\mbox{\rm
if}\;\;
\beta<\frac{3}{4}\;\;\mbox{\rm  and } (C(3));\\
 (y(1-y))^{2\gamma-1+\mu}, &&\mbox{\rm if}\;\;
\beta<\frac{3}{4},\;\;\mbox{\rm  and not}\;\; (C(3));\\
(y(1-y))^{2+2\gamma+\mu}, &&\mbox{\rm if}\;\; \beta\ge\frac{3}{4}.
\end{array}\right.
$$
We have the following result.
\begin{thm}\label{th-BK-general}
Under the conditions of {\rm Theorem} {\rm\ref{thm-HoHsing}} we have
with some $C_0>0$, as $n\to\infty$,
\begin{eqnarray}\label{main.result.2}
\lefteqn{\sup_{y\in (C_0\delta_n,1-C_0\delta_n)}\psi_3(y)\left|
R_n(y)-n^{-1}\sigma_{n,1}^{-1}\frac{f^{(1)}(Q(y))}{2}\left(\sum_{i=1}^nX_i\right)^2\right|}\\
&=&\left\{
\begin{array}{ll}
O_{\assubs}(d_{n,2}), &\mbox{\rm if}\;\; \beta\ge 2/3,\\
O_{a.s.}(c_n), &\mbox{\rm if}\;\; \beta<2/3  .
\end{array}
\right.\nonumber .\hspace*{6cm}
\end{eqnarray}
If $\gamma=1$ then the above estimate is valid on $(0,1)$.
\end{thm}

The (weighted) almost sure behavior of $R_n(\cdot)$ and (weighted)
convergence can be obtained in the same way as that of
$\tilde{R}_n(\cdot)$ in Corollaries \ref{cor.lil.uniform.BK} and
\ref{cor-weak-convergence}.
\subsection{Weak behavior of the general quantile process and its
consequences}\label{sec-weak-behaviour} Ho and Hsing's result
(\ref{quantiles-HoHsing}) would suggest that it should be possible
to approximate $q_n(y)$ at least on the expanding intevals,
$(n^{-1},1-n^{-1})$. However, as we will explain below, this is not
the case.

Let $\psi_4(y)=1$ or $y(1-y)$ according to $\beta<\frac{3}{4}$ or
$\beta\ge \frac{3}{4}$, respectively.
\begin{prop}\label{th-approx-unifquantile-generalquantile}
Assume {\rm(CsR1)-(CsR4)}. Then
\begin{equation}\label{quantiles-general-expo}
\sup_{y\in
(0,1)}\psi_4(y)|f(Q(y))q_n(y)-u_n(y)|=O_{a}(\sigma_{n,1}n^{-1}\ell
(n)),
\end{equation}
where $O_{a}=O_{\assubs}$ if $\gamma=1$, and $O_{a}=O_P$ if
$\gamma>1$.
\end{prop}
\begin{cor}\label{th-general-quantiles-weak-approx}
Assume {\rm(CsR1)-(CsR4)}. Then, under the conditions of either {\rm
Theorem} {\rm\ref{cor-reduction-quantiles}} or
{\rm\ref{th-reduction-quantiles-general}}, as $n\to\infty$,
$$
\sup_{y\in
(0,1)}\psi_1(y)f(Q(y))\left|q_n(y)+\sigma_{n,1}^{-1}\sum_{i=1}^nX_i\right|=O_P(n^{-(\beta-\frac{1}{2})}\ell
(n)).
$$
\end{cor}
\begin{cor}Assume {\rm(CsR1)-(CsR4)}. Then, under the conditions of either {\rm Theorem}
{\rm\ref{cor-reduction-quantiles}} or
{\rm\ref{th-reduction-quantiles-general}}, as $n\to\infty$,
$$
\sup_{y\in
(n^{-1},1-n^{-1})}(y(1-y))^{\nu}\left|q_n(y)+\sigma_{n,1}^{-1}\sum_{i=1}^nX_i\right|=o_P(1),
$$
where $$
\begin{array}{ll}
\nu>\gamma-(\beta-\frac{1}{2}), & \mbox{\rm if}\;\;
\beta<\frac{3}{4} \;\;\mbox{\rm and either (A(2)) or (C(2))};\\
\nu>2\gamma-\beta, &
\mbox{\rm if}\;\; \beta<\frac{3}{4}\;\; \mbox{\rm and neither (A(2)) nor (C(2))};\\
\nu>2\gamma-(\beta-\frac{1}{2}), & \mbox{\rm if}\;\;
\beta\ge\frac{3}{4} .
\end{array}
$$
\end{cor}

\indent From this result one obtains the following simultaneous
confidence bounds, which cover all the data available for $y\in
(n^{-1},1-n^{-1})$,
$$
Q_n(y)-\sigma_{n,1}n^{-1}c_{\nu}z_{\alpha}{(y(1-y))^{-\nu}}\le
Q(y)\le
Q_n(y)+\sigma_{n,1}n^{-1}c_{\nu}z_{\alpha}{(y(1-y))^{-\nu}},
$$
where $z_{\alpha}$ is the $(1-\alpha/2)$-quantile of the standard
normal law, and $$c_{\nu}=\sup_{y\in (0,1)}(y(1-y))^{\nu}.$$

Another consequence of Corollary
\ref{th-general-quantiles-weak-approx} is that for some
$k_n=k_n(\gamma,\beta)\to 0$, as $n\to\infty$,
$$
\sup_{y\in(k_n,1-k_n)}\left|q_n(y)+\sigma_{n,1}^{-1}\sum_{i=1}^nX_i\right|=o_P(1),
$$
and thus
\begin{equation}\label{weak-convergence-quantiles}
q_n(y)1_{\{y\in (k_n,1-k_n)\}}\convweak Z .
\end{equation}

\indent Optimally, one would hope to obtain weak convergence on
$(n^{-1},1-n^{-1})$, but this is not a good way to treat quantiles
in the LRD case at all. To see this, recall the subordinated
Gaussian model $Y_n=G(X_n)$. Take $G=F^{-1}\Phi$. For the uniform
sample quantile process $u_n(y)$ associated with the sequence
$\{Y_n,n\ge 1\}$ one obtains in the spirit of \cite[Proposition
2.2]{CSW2006} (see \cite{CsorgoKulik2006a} for a correct proof)
\begin{equation}\label{eq-reduction-subordination}
\sup_{y\in
(0,1)}\left|u_n(y)+\sigma_{n,1}^{-1}\phi(\Phi^{-1}(y))\sum_{i=1}^nX_i\right|=O_{P}(n^{-(\beta-\frac{1}{2})}\ell
(n)). \end{equation} Moreover, from \cite[Proposition 4.2]{CSW2006},
if the distribution $F$ of $Y=G(X)$ fulfills (CsR1)-(CsR3), then for some $k_n'\to 0$,
\begin{equation}\label{eq-approx-gaussian}
\sup_{y\in
(k_n',1-k_n')}\left|f(Q(y))q_n(y)+\sigma_{n,1}^{-1}\phi(\Phi^{-1}(y))\sum_{i=1}^nX_i\right|=O_P(n^{-(\beta-\frac{1}{2})}\ell
(n)),
\end{equation}
where $q_n(y)$ is the general quantile process associated with
$Y_n$. Thus,
\begin{equation}\label{weak-convergence-quantiles-1}
q_n(y)1_{\{y\in (k_n',1-k_n')\}}\convweak
\frac{\phi(\Phi^{-1}(y))}{f(Q(y))}Z ,
\end{equation}
provided $\frac{\phi(\Phi^{-1}(y))}{f(Q(y))}$ is uniformly bounded.
In particular, if $f$ is exponential, then this is not the case.
Consequently, we may have two LRD models, both with the same
covariance structure, both with the same exponential marginals, say,
so that in case of (\ref{model}) the general quantile process
converges, while in the subordinated Gaussian case it does not
converge (cf. (\ref{weak-convergence-quantiles}) and
(\ref{weak-convergence-quantiles-1}), respectively). On the other
hand, in both cases, the empirical processes have normal limits
sclaed by a deterministc function. In other words, subordination can
completely change convergence properties of quantile processes, even
if the empirical processes behave in the same way in the
subordinated and non-subordinated cases. The weight function
$(y(1-y))^{\nu}$ solves this problem somehow.
\subsubsection{Trimmed means}
In the model (\ref{model}), assume that $X_i$ are symmetric. From
(\ref{weak-convergence-quantiles}) one easily obtains
$$
\sigma_{n,1}^{-1}\left|\sum_{i=[nk_n]}^{[n(1-k_n)]}X_i\right|=\left|\int_{k_n}^{1-k_n}q_n(y)dy\right|\convdistr
|Z|.
$$
On the other hand, since $\Exp X_1=0$,
$\left|\int_0^1q_n(y)dy\right|=\sigma_{n,1}^{-1}\left|\sum_{i=1}^nX_i\right|\convdistr
|Z|$. If $k_n<l_n\to 0$ then the result remains true by considering
weak convergence in (\ref{weak-convergence-quantiles}) on
$(l_n,1-l_n)$ and then arguing as in the case of $k_n$. Summarizing,
\begin{cor}
Assume {\rm(CsR1)-(CsR4)} and that $X_i$ are symmetric. Let $k_n\le
l_n\to 0$. Then, under the conditions of either {\rm Theorem}
{\rm\ref{cor-reduction-quantiles}} or
{\rm\ref{th-reduction-quantiles-general}},
\begin{equation}\label{eq-trimmed-means}
\sigma_{n,1}^{-1}\sum_{i=[nl_n]}^{[n(1-l_n)]}X_i\convdistr Z .
\end{equation}
\end{cor}
The result (\ref{eq-trimmed-means}) states essentially that,
whatever trimming we consider, the deleted part is negligible.

However, it should be mentioned that this approach to the trimmed
sums is not the optimal one. The problem is considered in more details
in \cite{Kulik2006a} and \cite{KulikOuldHaye} via studying integral functionals of the empirical process (see e.g. \cite{Trimmed} for the description of the method in the i.i.d case).
\subsection{Remarks}\label{sec-remarks}
We start with pointing out some phenomena which are exclusive for
LRD sequences.
\begin{remark}{\rm
As mentioned in the Introduction, it was observed explicitly in
\cite{CSW2006} and can be concluded from \cite{Wu2005} that the
uniform Bahadur-Kiefer process (in case of \cite{CSW2006}) and,
under appropriate conditions, the general Bahadur-Kiefer process
(\cite{Wu2005}) converge in $D([y_0,y_1])$ for a particualt choice
of the parameter $\beta$. From our results we conclude that both
processes converge weakly in $D([0,1])$ if $\beta<\frac{3}{4}$. This
is striking difference compared to the i.i.d. case, for in the
latter case these processes cannot converge weakly (cf.
\cite{Kiefer1967}, \cite{Kiefer1970}). Considering pointwise
convergence, in the i.i.d. case the uniform and the general
Bahadur-Kiefer processes converge to the same limit (cf.
\cite{CsorgoSzyszkowicz1998} for a review). Here, the pointwise
limits are different, on account of different weak limits.}
\end{remark}
\begin{remark}\label{rem-uniform-see-general}{\rm
Unlike in the i.i.d case, to study the distance between the uniform
empirical and the uniform quantile processes, we need to control the
general quantile process, which can be done via controlling the
quantile and density quantile functions associated with $X_i$. The
reason for this is that the uniform quantile process contains
information regarding the marginal behavior of random variables
$X_i$. This is visible from Theorems \ref{cor-reduction-quantiles}
and \ref{th-reduction-quantiles-general} - the uniform quantile
process depends on the density-quantile function $f(Q(y))$
associated with $X_1$. As can be seen in
(\ref{eq-reduction-subordination}), this remains true in the
subordinated case $Y_i=G(X_i)$ as well, namely the uniform quantile
process contains information about the marginals of $X_i$, not of
$Y_i$. This has a impact on the behavior of general quantiles, as
described in Section \ref{sec-weak-behaviour}. }
\end{remark}

\indent We continue with some technical remarks concerning
assumptions and results above.

\begin{remark}{\rm
We comment on the different rates in our theorems, according to
different choices of $\beta$.

If $p=1$ then $a_n=o(d_{n,1})$, if $p=2$ (so that $\beta<3/4$), then
$d_{n,2}=o(a_n)$, and then optimal rates are attained in Theorems
\ref{cor-reduction-quantiles} and
\ref{th-reduction-quantiles-general}. Taking higher order expansions
($p\ge 3$) does not improve rates and requires additional
restrictions on $\beta$ and condtions conditions on $F$, either
(A(p)) or (C(p)).

Likewise, if $p=1,2$, then $c_n=o(d_{n,p})$. If $p=3$
($\beta<\frac{2}{3}$), then $d_{n,3}=o(c_n)$. Then we can identify
(but not prove !!) optimal rates in Theorems \ref{th-uniform},
\ref{th-reduction-quantiles-weighted}. We conjecture, that the bound
in Theorem \ref{th-uniform} (at least for $\beta<\frac{2}{3}$) is
valid without the $(\log n)^{1/2}$ term due to the following
conjecture.
\begin{conj} For
any $p\ge 1$,
$$
\limsup_{n\to\infty}\sigma_{n,p}^{-1}(\log\log n)^{-p/2}Y_{n,p}\as
c(\beta,p),
$$
\end{conj}
where $c(\beta,p)$ is as in Corollary \ref{cor.lil.uniform.BK}.

Further, on comparing Theorem \ref{th-BK-general} with
(\ref{BahadurKiefer-Wu}) we can see that the method in \cite{Wu2005}
leads to better rates for $\beta$ close to 1. We loose some rates
for $\beta$ close to 1, since then the error in the reduction
principle dominates. On the other hand, Wu's method is unlikely to
work when one wants to deal with approximations on the whole
interval $(0,1)$, which was our main goal. In fact, in view of a
weighted law of the iterated logarithm (see Lemma
\ref{lem-weighted-LIL-quantile}), it is not likely that in the case
$\beta\ge \frac{3}{4}$ the estimates on $(0,1)$ can be obtained with
optimal rates, unless the rate $d_{n,p}$ is improved. }
\end{remark}

\begin{remark}{\rm
Wu in his paper \cite{Wu2003} has in fact some weaker conditions on
$F_{\epsilon}$, than those stated in Theorem \ref{thm-HoHsing}.
Also, here, we avoid the boundary case $(p+1)(2\beta-1)=1$.
Furthermore, under stronger regularity conditions on the
distribution of $\epsilon_1$, the reduction principle (with worse
rates) for the empirical process remains true provided $\Exp
|\epsilon|^{2+\delta}<\infty$, $\delta>0$ (see
\cite{GiraitisSurgailis1999}). Thus, some of the results here remain
valid under the Giraitis and Surgailis conditions in
\cite{GiraitisSurgailis1999}. However, to prove Theorems
\ref{th-reduction-quantiles-general} and
\ref{th-reduction-quantiles-weighted} we require Lemma
\ref{lem-weighted-LIL} below, where the rates in the reduction
principle fo Theorem \ref{thm-HoHsing} are crucial. }
\end{remark}
\begin{remark}{\rm
We comment on assumptions (A(p)), (B) and (C(p)) on the distribution
function $F$. Note that $-(f\circ Q)^{(1)}(y)=J(y)$ is the so-called
score function (cf. e.g. \cite[p. 7]{CsorgoLN}), thus (A(1))
requires uniform boundness of the latter. This is not valid if one
takes the standard normal distribution for example. The assumptions
(A(p)), $p\ge 1$ are fulfilled if one takes the exponential,
logistic, or Pareto distribution $f(x)=\alpha( x^{1+\alpha})^{-1}$,
$x>1$, $\alpha>0$. Assumption (B) is fulfilled if one takes
exponential, logistic, or Pareto with $\alpha>1$. The latter
constrain $\alpha>1$ is relevant, since in view of Theorem
\ref{thm-HoHsing} we work under the condition $\Exp
\epsilon^4<\infty$ and, consequently, $\Exp X^4<\infty$. Further,
$(C(p))$, $p\ge 1$, is fulfilled in the Pareto case and for the
standard normal case. Thus, essentially, most of the "practical"
parametric families fulfill either (A(p)) or (C(p)).

Further, in the LRD case (\ref{model}) it is very unlikely that $f$
has bounded support (from either side). Moreover, to use of Theorem
\ref{thm-HoHsing}, we need $\Exp\epsilon=0$ and
$f_{\epsilon}=F^{'}_{\epsilon}$ to be smooth. Consequently, the same
properties are transferred to $X$ and its density $f$. Therefore, to
make use Theorem \ref{thm-HoHsing} and assumptions (A(p)) and (B)
simultaneously, we should consider the above comments for double
exponential or symmetric Pareto, appropriately smoothed around the
origin. Nevertheless, the main issue of assumptions (A(p)), (B) and
(C(p)) is the tail behavior. }
\end{remark}
\begin{remark}{\rm
As for the general quantile and the general Bahadur-Kiefer
processes, in order to obtain their approximations on the whole
interval, we assumed the monotonicity property (CsR4). In principe,
as in the i.i.d case, (cf. \cite{CsorgoCsorgoHorvathRevesz}), it
should be possible to obtain their approximations on the "practical"
interval $(n^{-1},1-n^{-1})$ without (CsR4).  }
\end{remark}
\begin{remark}{\rm
We now discuss the weights which appear in our theorems. As
mentioned in Remark \ref{rem-uniform-see-general}, the LRD sequences
based uniform quantile process "feels" the general quantile
function. In the i.i.d. case one knows that for $\mu>0$
$$
\limsup_{n\to\infty}\sup_{y\in (0,1)}(y(1-y))^{\mu}|Q(y)-Q_n^{\rm
iid}(y)|<\infty
$$
almost surely if and only if
$\int_{-\infty}^{\infty}|u|^{1/\mu}dF(u)<\infty$ (see \cite[p.
98]{CsorgoLN} for a tribute to David Mason in this regard).
Therefore, our weight functions $(y(1-y))^{\kappa}$, with some
$\kappa>0$, appear to be natural to use.

We also note that instead of the weight $(y(1-y))^{1+\kappa}$,
$\kappa>0$, we may consider $f^{\kappa'}(Q(y))$ as a weight
function, where $\kappa'$ depends on both $\kappa$ and $\gamma$. }
\end{remark}
\begin{remark}{\rm
In Theorem \ref{th-BK-general}, in case $\gamma>1$, the
approximation {\it in probability} remains valid on $(0,1)$ (see
also Proposition \ref{th-approx-unifquantile-generalquantile}). We
are not able to do this almost surely, since we do not have a
precise knowledge about the LRD behavior of order statistics (see
the proof of Proposition
\ref{th-approx-unifquantile-generalquantile}). }
\end{remark}
\begin{remark}\label{rem-full-picture}{\rm
The bound in Theorem \ref{cor-reduction-quantiles} is determined by
the behavior of the Bahadur-Kiefer process $\tilde R_n(y)$ (compare
Theorem \ref{cor-reduction-quantiles} with (\ref{exact.rate})). This
is somehow similar to the i.i.d. case. One knows that on an
appropriate probability space, $\sup_{y\in
(0,1)}|\alpha_n^{\mbox{\rm iid}}(y)-B_n(y)|=O_{\assubs}(n^{-1/2}\log
n)$, where $B_n(\cdot)$ are appropriate Brownian bridges. Further,
via (\ref{Kiefer}) we can see that with the same Brownian bridges we
have $\sup_{y\in (0,1)}|u_n^{\mbox{\rm
iid}}(y)-B_n(y)|=O_{\assubs}(n^{-1/4}(\log n)^{1/2}(\log\log
n)^{1/4})$. We may for example refer to \cite{CsorgoReveszbook} and
\cite{CsorgoSzyszkowicz1998} for more details.}
\end{remark}
\begin{remark}\label{rem-gaussian}{\rm
Recall, from Section \ref{sec-weak-behaviour}, our lines the
subordinated Gaussian case $Y=G(X)$. We have
$J_1(y)=-\phi(\Phi^{-1}(y))$, where $\phi$, $\Phi$ are the standard
normal density and distribution function. Cs\"{o}rg\H{o},
Szyszkowicz and Wang in \cite{CSW2006} proved their Proposition 2.2
assuming (cf. also their Remark 2.1) their Assumption A. However,
what is really used in their proof is that $J_1$ has, in particular,
uniformly bounded first order derivative, which is not true, since
$J_1'(y)=-\Phi^{-1}(y)$. Consequently, their Proposition 2.2 and all
its consequences in their Sections 2.1 and 2.2 are valid only if one
restricts them to intervals $[y_0,y_1]$, or assumes that $Y=G(X)$
has finite support. This actually is the reason that we considered
assumptions (A(p)), (B) and/or weighted approximations. Clearly, the
non-subordinated Gaussian case can be treated as in the setting of
Theorems \ref{th-reduction-quantiles-general},
\ref{th-reduction-quantiles-weighted} and \ref{th-BK-general} with
$\gamma=1$ (recall that (C(p)) holds in the Gaussian case). For the
general treatment we refer to \cite{CsorgoKulik2006a}.

Also, as noted already in our Section \ref{sec-general-BK}, results
for the general Bahadur-Kiefer process cannot be concluded from an
approximation of the latter by the uniform one. Hence, the proposed
proofs for Theorems 4.1, 4.2 of \cite{CSW2006} via the invariance
principle of Proposition 4.2 cannot work and, in view of
\cite{Wu2005}, the claimed limiting processes can at best be correct
if multiplied by 1/2.

In Section 3 of \cite{CSW2006} the authors consider
$V_n(t)=2\sigma_{n,1}^{-1}n\int_0^t\tilde R_n(y)dy$ and
$Q_n(t)=V_n(t)-\alpha_n^2(t)$, the so-called uniform Vervaat and
Vervaat Error processes. As a consequence of our comments so far on
paper \cite{CSW2006}, we note that the results in this section are
valid only if $G(X)$ has finite support. An extension is possible if
one has assumptions like (A(p)) and (B). This, however, is out of
the scope of this paper. }
\end{remark}
\section{Proofs}\label{sec-proofs}
\subsection{Preliminary results}
We recall the following law of the iterated logarithm for partial
sums $\sum_{i=1}^nX_i$ (see, e.g., \cite{WangLinGulati}):
\begin{equation}\label{step.5}
\limsup_{n\to\infty}\sigma_{n,1}^{-1}(\log\log
n)^{-1/2}\left|\sum_{i=1}^nX_i\right|\as c(\beta,1),
\end{equation}
where $c(\beta,1)$ is defined in Corollary \ref{cor.lil.uniform.BK}.
\begin{lem}\label{lem.strong.bound}
Let $p\ge 1$ be an arbitrary integer such that $p<(2\beta-1)^{-1}$.
Then, as $n\to\infty$,
\begin{equation}\label{conjecture}
Y_{n,p}=O_{\assubs}(\sigma_{n,p}(\log n)^{1/2}\log\log n).
\end{equation}
\end{lem}
{\it Proof.} Let $B_n^2=\sigma_{n,p}^2\log n (\log\log n)^2$. By
(\ref{eq-varance-behaviour}), \cite[Lemma 4]{Wu2005} and Karamata's
Theorem we have
\begin{eqnarray*}
\lefteqn{\hspace*{-2cm}
\left\|\frac{|Y_{n,p}|}{B_{2^d}}\right\|_2^2\le
\frac{1}{B_{2^d}}\left(\sum_{j=0}^{d}2^{(d-j)/2}\sigma_{2^j,p}\right)^2\le
\frac{2^d}{B_{2^d}}\left(\sum_{j=0}^d2^{j(1-p(2\beta-1))/2}L_0^p(2^j)\right)^2}\\
&\sim & \frac{2^d}{B_{2^d}}2^{2d-dp(2\beta-1)}L_0^{2p}(2^d)\sim
d^{-1}(\log d)^{-2}.
\end{eqnarray*}
Therefore, the result follows by the Borel-Cantelli lemma.
 \koniec

As an easy consequence of (\ref{step.5}) and (\ref{conjecture}) we
obtain the next result.
\begin{lem}
Let $p\ge 1$ be an arbitrary integer such that $p<(2\beta-1)^{-1}$.
We have
\begin{equation}\label{step.3}
\limsup_{n\to\infty}\sigma_{n,1}^{-1}(\log\log n)^{-1/2}\sup_{y\in
(0,1)}|\tilde V_{n,p}(y)|\as c(\beta,1).
\end{equation}
\end{lem}

\indent Using Theorem \ref{thm-HoHsing} and the same argument as in
the proof of Lemma \ref{lem.strong.bound}, we obtain
\begin{eqnarray*}
\lefteqn{\sigma_{n,p}^{-1}\sup_{x\in
{\matR}}|S_{n,p}(x)|}\\
&=&\left\{\begin{array}{ll} O_{\assubs}(n^{-(\frac{1}{2}-p(\beta-\frac{1}{2}))}L_0^{-p}(n)(\log n)^{5/2}(\log\log n)^{3/4}), & (p+1)(2\beta-1)>1\\
O_{\assubs}(n^{-(\beta-\frac{1}{2})}L_0(n)(\log n)^{1/2}(\log \log
n)^{3/4}), & (p+1)(2\beta-1)<1
\end{array}\right. .
\end{eqnarray*}
Since (see (\ref{eq-varance-behaviour}))
\begin{equation}\label{variances}
\frac{\sigma_{n,p}}{\sigma_{n,1}}\sim
n^{-(\beta-\frac{1}{2})(p-1)}L_0^{p-1}(n),
\end{equation}
we obtain
\begin{eqnarray*}
\lefteqn{\sup_{x\in {\matR}}|\beta_n(x)+\sigma_{n,1}^{-1}V_{n,p}(x)|=}\\
& = & \frac{\sigma_{n,p}}{\sigma_{n,1}}\sup_{x\in
{\matR}}\left|\sigma_{n,p}^{-1}\sum_{i=1}^n(1_{\{X_i\le
x\}}-F(x))+\sigma_{n,p}^{-1}V_{n,p}(x)\right|=o_{\assubs}(d_{n,p}).
\end{eqnarray*}
Consequently, via $\{\alpha_n(y),y\in (0,1)\}=\{\beta_n(Q(y)),y\in
(0,1)\}$,
\begin{equation}\label{approx-unif-empirical}
\sup_{y\in (0,1)}|\alpha_n(y)+\sigma_{n,1}^{-1}\tilde
V_{n,p}(y)|=O_{\assubs}(d_{n,p}).
\end{equation}
\begin{remark}\label{rem-constants}{\rm
For convenient reference, we collect here various relations between
constants. Recall that $d_{n,2}=o(a_n)$ provided
$\beta<\frac{3}{4}$, and $d_{n,3}=o(c_n)$, provided
$\beta<\frac{2}{3}$. Further, $\sigma_{n,1}^{-1}b_{n,p}=o(d_{n,p})$.
It is not necessarily true that $\sigma_{n,1}^{-1}=o(d_{n,p})$, but
it is always true that $\sigma_{n,1}^{-1}=o(a_n)$. }
\end{remark}
\subsection{Proof of Theorems \ref{cor-reduction-quantiles} and \ref{th-uniform}}
First, we bound the distance between the uniform empirical and
uniform quantile processes.
\begin{lem}\label{lem-distance-unemp-unquan}
Let $p\ge 1$ be an arbitrary integer such that $p<(2\beta-1)^{-1}$.
Assume {\rm(A(p))}. Under the conditions of {\rm Theorem
{\rm\ref{thm-HoHsing}}} we have, as $n\to\infty$,
$$
\sup_{y\in
(0,1)}|u_n(y)-\alpha_n(y)|=O_{\assubs}(a_n)+O_{\assubs}(d_{n,p}).
$$
\end{lem}
{\it Proof.} Note that
\begin{eqnarray}
u_n(y)&=&\sigma_{n,1}^{-1}n(E_n(U_{n}(y))-U_{n}(y))-\sigma_{n,1}^{-1}n(E_n(U_{n}(y))-y)\\
&= &
\sigma_{n,1}^{-1}n(E_n(U_{n}(y))-U_{n}(y))+O_{\assubs}(\sigma_{n,1}^{-1})=\alpha_n(U_n(y))+O(\sigma_{n,1}^{-1})
.\label{geometric-property}\nonumber
\end{eqnarray}
Thus, by (\ref{approx-unif-empirical}),
\begin{eqnarray}
\lefteqn{\sup_{y\in
(0,1)}|u_n(y)-\alpha_n(y)|}\label{step.8}\\
&= & \sup_{y\in
(0,1)}|\alpha_n(U_{n}(y))-\alpha_n(y)|+O_{\assubs}(\sigma_{n,1}^{-1})\nonumber\\
&\le & \sigma_{n,1}^{-1}\sup_{y\in (0,1)}|\tilde V_{n,p}(y)-\tilde
V_{n,p}(U_{n}(y))|+O_{\assubs}(\sigma_{n,1}^{-1})+O_{\assubs}(d_{n,p})\nonumber
.
\end{eqnarray}
Accordingly, in view of Assumptions (A(p)), (B), we have to control
\begin{equation}\label{step.1}
\sup_{y\in
(0,1)}\left|f(Q(y))-f(Q(U_{n}(y)))\right|\left|\sum_{i=1}^{n}X_i\right|\le
C\sup_{y\in (0,1)}|y-U_n(y)|\left|\sum_{i=1}^{n}X_i\right|
\end{equation}
and
\begin{eqnarray}\label{step.2}
\lefteqn{\sup_{y\in
(0,1)}\sum_{r=2}^p\left|f^{(r-1)}(Q(y))-f^{(r-1)}(Q(U_{n}(y)))\right||Y_{n,r}|}\\
&\le& C\sup_{y\in
(0,1)}|y-U_n(y)|\left|\sum_{r=2}^pY_{n,r}\right|.\hspace*{5cm}\nonumber
\end{eqnarray}
From (\ref{step.3}) and (\ref{approx-unif-empirical}) one
obtains
$$
\limsup_{n\to\infty}(\log\log n)^{1/2}\sup_{y\in
(0,1)}|\alpha_n(y)|\as c(\beta,1).
$$
Consequently, as $n\to\infty$,
\begin{eqnarray}
\sup_{y\in (0,1)}|y-U_n(y)|&=&\sup_{y\in
(0,1)}\sigma_{n,1}n^{-1}|u_n(y)|=\sup_{y\in
(0,1)}\sigma_{n,1}n^{-1}|\alpha_n(y)|\nonumber\\
&=&O_{\assubs}(\sigma_{n,1}n^{-1}(\log\log
n)^{1/2})=O_{\assubs}(a_n). \label{step.4}
\end{eqnarray}
Therefore, on combining (\ref{step.5}), (\ref{step.1}),
(\ref{step.4}), as $n\to\infty$, one obtains
\begin{equation}\label{step.6}
\sup_{y\in
(0,1)}\sigma_{n,1}^{-1}|f(Q(y))-f(Q(U_n(y)))|\left|\sum_{i=1}^nX_i\right|=O_{\assubs}(a_n).
\end{equation}
Having (\ref{conjecture}), (\ref{step.2}) and (\ref{step.6}), as
$n\to\infty$, we conclude
\begin{equation}\label{step.7}
\sup_{y\in (0,1)}\sigma_{n,1}^{-1}|\tilde V_{n,p}(y)-\tilde
V_{n,p}(U_n(y))|=O_{\assubs}(a_n).
\end{equation}
Thus, by (\ref{step.8}) and (\ref{step.7}), as $n\to\infty$,
\begin{eqnarray*}
\sup_{y\in
(0,1)}|u_n(y)-\alpha_n(y)|=O_{\assubs}(a_n)+O(\sigma_{n,1}^{-1})+O_{\assubs}(d_{n,p}),
\end{eqnarray*}
and hence the result follows. \koniec

If $\beta\ge 3/4$, take $p=1$ and assume (A(1)). If $\beta<3/4$,
take $p=2$ and assume (A(2)). As a consequence of Lemma
\ref{lem-distance-unemp-unquan}, (\ref{conjecture}),
(\ref{approx-unif-empirical}) and Remark \ref{rem-constants} we
obtain (\ref{eq-uniform-quantile-approximation}).
\subsubsection{Proof of Theorem \ref{th-uniform}}
 In Lemma
\ref{lem-distance-unemp-unquan} we have a bound on the distance
between the uniform empirical and the uniform quantile processes,
but it does not say anything about its optimality. To obtain this
note, that for any $1\le p<(2\beta-1)^{-1}$ we have by
(\ref{approx-unif-empirical}) and as in (\ref{geometric-property})
\begin{eqnarray}
\lefteqn{\sup_{y\in
(0,1)}|\alpha_n(y)-u_n(y)+\sigma_{n,1}^{-1}(\tilde
V_{n,p}(y)-\tilde{V}_{n,p}(U_n(y)))| }\nonumber\\
&\le & \sup_{y\in
(0,1)}|\alpha_n(y)-\alpha_n(U_n(y))+\sigma_{n,1}^{-1}(\tilde
V_{n,p}(y)-\tilde{V}_{n,p}(U_n(y)))|
\label{BahadurKiefer-approximation-step}\\&&+\sup_{y\in
(0,1)}|\alpha_n(U_n(y))-u_n(y)|=O_{\assubs}(d_{n,p})+O_{\assubs}(\sigma_{n,1}^{-1})\nonumber.
\end{eqnarray}
Now, it is sufficient to deal with the process $(\tilde
V_{n,p}(y)-\tilde{V}_{n,p}(U_n(y)))$. We approximate this process
via several lemmas.
\begin{lem}\label{lem.1}
Let $p\ge 1$ be an arbitrary integer such that $p<(2\beta-1)^{-1}$.
Assume {\rm(A(p))} and {\rm(B)}. Under the conditions of {\rm
Theorem {\rm\ref{thm-HoHsing}}} we have as $n\to\infty$,
$$
\sup_{y\in (0,1)}\left|\tilde V_{n,1}(y)-\tilde
V_{n,1}(U_n(y))+\frac{f^{(1)}(Q(y))}{f(Q(y))}\frac{\tilde{V}_{n,p}(y)}{n}\sum_{i=1}^nX_i\right|=O_{\assubs}(b_n)+O_{\assubs}(b_{n,p}).
$$
\end{lem}
{\it Proof.} Applying second order Taylor expansion and recalling
that $(f\circ Q)^{(1)}(y)=\frac{f^{(1)}(Q(y))}{f(Q(y))}$, one
obtains
\begin{eqnarray*}
\lefteqn{\sup_{y\in
(0,1)}\left|(f(Q(y))-f(Q(U_n(y))))\sum_{i=1}^nX_i+n^{-1}\frac{f^{(1)}(Q(y))}{f(Q(y))}\tilde
V_{n,p}(y)\sum_{i=1}^nX_i\right|}\\
& \le & \sup_{y\in
(0,1)}\left|\frac{f^{(1)}(Q(y))}{f(Q(y))}\sigma_{n,1}n^{-1}\sum_{i=1}^nX_i\left(u_n(y)+\sigma_{n,1}^{-1}\tilde
V_{n,p}(y)\right)\right|\\
&&+ \sup_{y\in (0,1)}|(f\circ Q)^{(2)}(y)|\sup_{y\in
(0,1)}(y-U_n(y))^2\left|\sum_{i=1}^nX_i\right|
\\
& = & O_{\assubs}(\sigma_{n,1}n^{-1}a_n\sigma_{n,1}(\log\log
n)^{1/2})+O_{\assubs}(\sigma_{n,1}n^{-1}d_{n,p}\sigma_{n,1}(\log\log
n)^{1/2})\\
&&+O_{\assubs}(\sigma_{n,1}^3n^{-2}(\log\log
n)^{3/2})=O_{\assubs}(b_n)+O_{\assubs}(b_{n,p}).
\end{eqnarray*}
The above bound follows from (\ref{step.5}),
(\ref{approx-unif-empirical}), (\ref{step.4}) and
(\ref{eq-uniform-quantile-approximation}) Theorem
\ref{cor-reduction-quantiles}.
 \koniec
\begin{lem}\label{lem.3}
Let $p\ge 1$ be an arbitrary integer such that $p<(2\beta-1)^{-1}$.
Assume {\rm(A(p))} and {\rm(B)}. Under the conditions of {\rm
Theorem {\rm\ref{thm-HoHsing}}} we have as $n\to\infty$,
$$
\sup_{y\in (0,1)}n^{-1}\frac{f^{(1)}(Q(y))}{f(Q(y))}|\tilde
V_{n,p}(y)-\tilde
V_{n,1}(y)|\left|\sum_{i=1}^nX_i\right|=O_{\assubs}(b_n(\log
n)^{1/2}).
$$
\end{lem}
{\it Proof.} We have
\begin{eqnarray*}
\lefteqn{\sup_{y\in (0,1)}n^{-1}\left|\tilde V_{n,p}(y)-\tilde
V_{n,1}(y)\right|\left|\sum_{i=1}^nX_i\right|}\\
&\le & \sup_{y\in
(0,1)}|f^{(1)}(Q(y))|n^{-1}|Y_{n,2}|\left|\sum_{i=1}^nX_i\right|+O_{\assubs}\left(n^{-1}\left|\sum_{r=3}^pY_{n,r}\right|\left|\sum_{i=1}^nX_i\right|\right).
\end{eqnarray*}
Using (\ref{step.5}), (\ref{conjecture}), we obtain the result.
\koniec

Similarly to Lemma \ref{lem.3}, the next result holds true as well.
\begin{lem}\label{lem.2}
Let $p\ge 1$ be an arbitrary integer such that $p<(2\beta-1)^{-1}$.
Assume {\rm(A(p))} and {\rm(B)}. Under the conditions of {\rm
Theorem {\rm\ref{thm-HoHsing}}} we have as $n\to\infty$,
$$
\sup_{y\in (0,1)}|\tilde V_{n,1}(y)-\tilde V_{n,1}(U_n(y))-(\tilde
V_{n,p}(y)-\tilde V_{n,p}(U_n(y)))|=O_{\assubs}(b_n(\log n)^{1/2}).
$$
\end{lem}

From Lemmas \ref{lem.1}, \ref{lem.3}, \ref{lem.2} we obtain
\begin{cor}\label{cor.1}
Let $p\ge 1$ be an arbitrary integer such that $p<(2\beta-1)^{-1}$.
Assume {\rm(A(p))} and {\rm(B)}. Under the conditions of {\rm
Theorem {\rm\ref{thm-HoHsing}}} we have as $n\to\infty$,
\begin{eqnarray*}
\lefteqn{\hspace*{-3cm}\sup_{y\in (0,1)}\left|\tilde
V_{n,p}(y)-\tilde
V_{n,p}(U_n(y))+n^{-1}\frac{f^{(1)}(Q(y))}{f(Q(y))}\tilde
V_{n,1}(y)\sum_{i=1}^nX_i\right|}\\
&= &O_{\assubs}(b_n(\log n)^{1/2})+O_{as}(b_{n,p}).
\end{eqnarray*}
\end{cor}
\indent Recall that $\tilde{R}_n(y)=\alpha_n(y)-u_n(y)$. Then, by
(\ref{BahadurKiefer-approximation-step}),
$$
\sup_{y\in (0,1)}|\tilde R_n(y)+\sigma_{n,1}^{-1}(\tilde
V_{n,p}(y)-\tilde
V_{n,p}(U_n(y)))|=O_{\assubs}(d_{n,p})+O_{\assubs}(\sigma_{n,1}^{-1}).
$$
Consequently, via Corollary \ref{cor.1}, \begin{eqnarray*}
\lefteqn{\sup_{y\in (0,1)}\left|\tilde
R_n(y)-n^{-1}\sigma_{n,1}^{-1}\frac{f^{(1)}(Q(y))}{f(Q(y))}\tilde
V_{n,1}(y)\sum_{i=1}^nX_i\right|}\\
& =& O_{\assubs}(d_{n,p})+O_{\assubs}(\sigma_{n,1}^{-1}b_n(\log
n)^{1/2})+O_{\assubs}(\sigma_{n,1}^{-1}b_{n,p})+O_{\assubs}(\sigma_{n,1}^{-1})\\
&=&O_{\assubs}(d_{n,p})+O_{\assubs}(c_n)+O_{\assubs}(\sigma_{n,1}^{-1}).
\end{eqnarray*}
If $\beta\ge 2/3$, then the bound is $O_{\assubs}(d_{n,2})$ on
assuming (A(2)). If $\beta<2/3$, taking $p=3$, via Remark
\ref{rem-constants}, we obtain the statement
(\ref{quantile-approx-weighted}) of Theorem \ref{th-uniform}.
\koniec
\subsection{Proof of the optimality in Theorem \ref{cor-reduction-quantiles}}
If $\beta<\frac{3}{4}$, then the dominating term in Theorem
\ref{cor-reduction-quantiles} is $O_{\assubs}(a_n)$.

Fix $y=y_0$. Via (\ref{approx-unif-empirical}) and as in
(\ref{exact.rate}) we obtain
\begin{eqnarray*}
\lefteqn{\limsup_{n\to\infty}\sigma_{n,1}^{-1}n(\log\log
n)^{-1}\left|u_n(y_0)+\sigma_{n,1}^{-1}f(y_0)\sum_{i=1}^nX_i\right|}\\
&= & \limsup_{n\to\infty}n(\sigma_{n,1}\log\log
n)^{-1}|u_n(y_0)-\alpha_n(y_0)+(\alpha_n(y_0)+\sigma_{n,1}^{-1}\tilde
V_{n,p}(y_0))|\\
&= & c(\beta,1)|f^{(1)}(Q(y_0))|.
\end{eqnarray*}
Therefore, via (\ref{eq-uniform-quantile-approximation}), for any
$y_0\in (0,1)$,
$$
c(\beta,1)|f^{(1)}(Q(y_0))|\le\limsup_{n\to\infty}\frac{n}{\sigma_{n,1}\log\log
n}\sup_{y\in (0,1)}\left|u_n(y)+\frac{\tilde
V_{n,p}(y)}{\sigma_{n,1}}\right|= O_{\assubs}(1)
$$
which means that the bound is optimal. \koniec
\subsection{Proof of Theorems \ref{th-reduction-quantiles-general} and \ref{th-reduction-quantiles-weighted} }\label{section-furtherproofs}
\subsubsection{Properties of the density-quantile function}
Note that under an appropriate smoothness of $f$, (CsR3) is
equivalent to
\begin{itemize}
\item[{\rm (CsR3(i))}] $f(Q(y))\sim y^{\gamma_1}L_1(y^{-1})$, as $y\downarrow
0$,
\item[{\rm (CsR3(ii))}] $f(Q(1-y))\sim (1-y)^{\gamma_2}L_2((1-y)^{-1})$, as $y\uparrow
1$,
\end{itemize}
for some numbers $\gamma_1,\gamma_2>0$ and some slowly varying
functions $L_1,L_2$. The parameter $\gamma$ in (CsR3) and
$\gamma_1,\gamma_2$ are related as $\gamma=\gamma_1\wedge \gamma_2$
(see \cite{CsorgoZitikis1999}). Let $\gamma_0=\gamma_1\vee\gamma_2$.
Under (CsR3(i)) and (CsR3(ii)) we have for any $\mu>0$,
\begin{equation}\label{eq-y-fQy}
\sup_{y\in (0,1)}\frac{(y(1-y))^{\gamma+\mu}}{f(Q(y))}=O(1).
\end{equation}
Further, note that if $0<\gamma_1<1$ ($0<\gamma_2<1$) then $F$ has
bounded support from the left (from the right) (see
\cite{Parzen1979}). Thus, we assume without loss of generality that
both $\gamma_1$ and $\gamma_2$ are not smaller than $1$. In this
case, for any $\varepsilon>0$,
\begin{equation}
\label{eq-fQy-y} f(Q(y))=O(y^{1-\varepsilon}),\quad y\to 0 .
\end{equation}
Note also, that (CsR3(i)) and (CsR3(ii)) together with $\gamma_0>1$
imply that for any $\mu>0$,
\begin{equation}\label{eq-one-denominator}
\sup_{y\in (0,1)}\frac{|f^{(1)}(Q(y))|}{f(Q(y))}(y(1-y))^{\mu}=O(1).
\end{equation}
Further, by \cite[p. 116]{Parzen1979},
\begin{equation}\label{eq-secondorderdensityquantile}
(f\circ Q)^{(2)}(y)\sim \kappa
\frac{(f^{(1)}(Q(y)))^2}{f^3(Q(y))}\end{equation} as $y\to 0$. The
parameter $\kappa$ is positive if $\gamma_1>1$ or $\kappa=0$ if
$\gamma_1=1$. A similar consideration applies to the upper tail.
\subsubsection{Weighted law of the iterated logarithm}
From (\ref{step.5}), (\ref{approx-unif-empirical}), (\ref{eq-fQy-y})
and $\delta_n^{-1/2}d_{n,p}=O(1)$ if $p\ge 2$ (i.e.
$\beta<\frac{3}{4}$) one obtains
\begin{lem}\label{lem-weighted-LIL}
Let $\beta<\frac{3}{4}$. Under the conditions of {\rm Theorem
{\rm\ref{thm-HoHsing}}}, as $n\to\infty$,
$$
\sup_{y\in
(\delta_n,1-\delta_n)}\frac{|\alpha_n(y)|}{(y(1-y))^{1/2}}=O_{\assubs}((\log\log
n)^{1/2}).
$$
\end{lem}
Using now the same argument as in \cite[Theorem
2]{CsorgoRevesz1978}, we obtain a corresponding result for the
linear LRD based uniform quantile process.
\begin{lem}\label{lem-weighted-LIL-quantile}
Let $\beta<\frac{3}{4}$. Under the conditions of {\rm Theorem
{\rm\ref{thm-HoHsing}}}, with some $C_0>0$, as $n\to\infty$,
$$
\sup_{y\in
(C_0\delta_n,1-C_0\delta_n)}\frac{|u_n(y)|}{(y(1-y))^{1/2}}=O_{\assubs}((\log\log
n)^{1/2}).
$$
\end{lem}

\indent From Lemma \ref{lem-weighted-LIL-quantile}, by the same
argument as in \cite[Theorem 3]{CsorgoRevesz1978}, as $n\to\infty$,
\begin{equation}\label{eq-small-increments-quantile}
\sup_{y\in (0,\delta_n)}|u_n(y)|=O_{\assubs}(a_n),
\end{equation}
provided $\beta<\frac{3}{4}$. Further, via (\ref{step.5}),
(\ref{approx-unif-empirical}) and (\ref{eq-fQy-y}), as $n\to\infty$,
we obtain for arbitrary $\beta\in (1/2,1)$ and $1\le
p<(2\beta-1)^{-1}$,
$$\sup_{y\in
(0,\delta_n)}|\alpha_n(y)|=O_{\assubs}(\delta_n^{1-\varepsilon}(\log\log
n)^{1/2})+O_{\assubs}(d_{n,p})=O_{\assubs}(a_n)+O_{\assubs}(d_{n,p}).$$

Recall (\ref{step.4}). Let $\theta=\theta_n(y)$ be such that
$|\theta-y|\le
\sigma_{n,1}n^{-1}|u_n(y)|=O_{\assubs}(n^{-(\beta-\frac{1}{2})}L_0(n)(\log\log
n)^{1/2})$. Arguing as in \cite[Theorem 3]{CsorgoRevesz1978},
uniformly for $y\in (C_0\delta_n,1-C_0\delta_n)$, as $n\to\infty$,
\begin{equation}\label{eq-bound-on-y-theta}
\frac{y(1-y)}{\theta(1-\theta)}=O_{\assubs}(1) .
\end{equation}
\subsubsection{Proof of Theorem
\ref{th-reduction-quantiles-general}} First, we need estimates which
will replace a part of the proof of Lemma
\ref{lem-distance-unemp-unquan}. All random variables $\theta$ below
are as in (\ref{eq-bound-on-y-theta}).
\begin{lem}\label{lem-extension-1}
Let $p\ge 1$ be an arbitrary integer such that $p<(2\beta-1)^{-1}$
and assume that {\rm(C(p))} is fulfilled. Under the conditions of
{\rm Theorem} {\rm\ref{th-reduction-quantiles-general}}, for any
$r=0,\ldots,p-1$, as $n\to\infty$,
$$
\sup_{y\in
(C_0\delta_n,1-C_0\delta_n)}\psi_1(y)|f^{(r)}(Q(y))-f^{(r)}(Q(U_n(y)))|=O_{\assubs}(n^{-(\beta-\frac{1}{2})}L_0(n)(\log\log
n)^{1/2}).
$$
\end{lem}
{\it Proof.} Let $\beta<\frac{3}{4}$. Take first
$\psi_1(y)=(y(1-y))^{\gamma-\frac{1}{2}+\mu}$. Taking a first order
Taylor expansion and bearing in mind that $f^{(r+1)}$ are uniformly
bounded, we have
$$
\psi_1(y)|f^{(r)}(Q(y))-f^{(r)}(Q(U_n(y)))|=\frac{(\theta(1-\theta))^{\gamma+\mu}}{f(Q(\theta))}\left(\frac{y(1-y)}{\theta(1-\theta)}\right)^{\gamma+\mu}\frac{|y-U_n(y)|}{(y(1-y))^{1/2}}.
$$
Further, under the condition (C(p)),
\begin{eqnarray*}
\lefteqn{
|f^{(r)}(Q(y))-f^{(r)}(Q(U_n(y)))|}\\
&=&\frac{f^{(r+1)}(Q(\theta))}{f(Q(\theta))}(\theta(1-\theta))^{1/2}\left(\frac{y(1-y)}{\theta(1-\theta)}\right)^{1/2}\frac{|y-U_n(y)|}{(y(1-y))^{1/2}}.
\end{eqnarray*}
Thus, the result follows by Lemma \ref{lem-weighted-LIL-quantile},
(\ref{eq-y-fQy}) and (\ref{eq-bound-on-y-theta}).

If $\beta\ge\frac{3}{4}$, assume (C(1)). We use the appropriate form
of $\psi_1$, (\ref{eq-y-fQy}) and (\ref{eq-bound-on-y-theta}).
\koniec

From Lemma \ref{lem-extension-1}, and exactly as in the proof of
Lemma \ref{lem-distance-unemp-unquan}, as $n\to\infty$,
$$
\sup_{y\in
(C_0\delta_n,1-C_0\delta_n)}\psi_1(y)|u_n(y)-\alpha_n(y)|=O_{\assubs}(a_n)+O_{\assubs}(d_{n,p}).
$$
Consequently, by (\ref{eq-small-increments-quantile}) and the
comment below it, as $n\to\infty$, we have for $\beta<\frac{3}{4}$
and $p<(2\beta-1)^{-1}$,
$$
\sup_{y\in
(0,1)}\psi_1(y)|u_n(y)-\alpha_n(y)|=O_{\assubs}(a_n)+O_{\assubs}(d_{n,p}).
$$
The same estimates are valid for $\beta\ge \frac{3}{4}$, since in
this case $\psi_1(y)=O(y)$. Consequently,
(\ref{quantile-approx-weighted}) follows.\koniec
\subsubsection{Proof of Theorem \ref{th-reduction-quantiles-weighted}}
First, we show that Lemma \ref{lem.1} remains valid when multiplying
by $\psi_2(y)$.

From (\ref{eq-one-denominator}), Theorem
\ref{th-reduction-quantiles-general} and estimating as in Lemma
\ref{lem.1}, as $n\to\infty$, we conclude
\begin{eqnarray}
\lefteqn{\hspace*{-4cm}\sup_{y\in
(0,1)}(y(1-y))\psi_1(y)\left|\frac{f^{(1)}(Q(y))}{f(Q(y))}\sigma_{n,1}n^{-1}\sum_{i=1}^nX_i\left(u_n(y)+\sigma_{n,1}^{-1}\tilde
V_{n,p}(y)\right)\right|}\nonumber\\
&=&O_{\assubs}(b_n)+O_{\assubs}(c_n)\label{eq-oszacowania-1}.
\end{eqnarray}
In view of (\ref{eq-secondorderdensityquantile}), for the term in
Lemma \ref{lem.1} involving $(f\circ Q)^{(2)}(y)$, we estimate
\begin{eqnarray}
\lefteqn{(y(1-y))^{\mu}\frac{(f^{(1)}(Q(\theta)))^2}{f^3(Q(\theta))}(y-U_n(y))^2\left|\sum_{i=1}^nX_i\right|}\nonumber\\
&=&\left(\frac{f^{(1)}(Q(\theta))}{f^2(Q(\theta))}\theta(1-\theta)\right)^2\frac{f(Q(\theta))}{(\theta(1-\theta))^{1-\mu}}\left(\frac{y(1-y)}{\theta(1-\theta)}\right)^{1+\mu}\frac{(y-U_n(y))^2}{y(1-y)}\left|\sum_{i=1}^nX_i\right|\nonumber\\
&= & O_{\assubs}(b_n)\label{eq-oszacowania},
\end{eqnarray}
uniformly for $y\in (C_0\delta_n,1-C_0\delta_n)$, on account of
(CsR3), (\ref{eq-fQy-y}), (\ref{eq-bound-on-y-theta}), Lemma
\ref{lem-weighted-LIL-quantile} and (\ref{step.5}). A similar
argument yields the same bound for the right tail.

Further, as $n\to\infty$,
\begin{eqnarray}
\lefteqn{\sup_{y\in (0,C_0\delta_n)}(y(1-y))^{1+\mu}|\tilde
V_{n,1}(y)-\tilde V_{n,1}(U_n(y))|}\label{eq-oszacowania-2}\\
&\le & C_0\delta_n^{1+\mu}\sup_{y\in (0,1)}f(Q(y))|\sum_{i=1}^nX_i|
= O_{\assubs}(\delta_n^{1+\mu}\sigma_{n,1}(\log\log
n)^{1/2})=o_{\assubs}(b_n)\nonumber.
\end{eqnarray}
and by (\ref{eq-one-denominator})
\begin{eqnarray}
\lefteqn{\sup_{y\in (0,C_0\delta_n)}(y(1-y))^{1+\mu}\left|\frac{f^{(1)}(Q(y))}{f(Q(y))}\right|\left|n^{-1}\tilde V_{n,p}\sum_{i=1}^nX_i\right|}\nonumber\\
&= & \delta_n^{1+\mu/2}\sup_{y\in
(0,C_0\delta_n)}(y(1-y))^{\mu/2}\left|\frac{f^{(1)}(Q(y))}{f(Q(y))}\right|O_{\assubs}\left(\left(\sum_{i=1}^nX_i\right)^2/n\right)\nonumber\\
&=&O_{\assubs}(\delta_n^{1+\mu/2}\sigma_{n,1}^2n^{-1}\log\log
n)=O_{\assubs}(b_n)\label{eq-oszacowania-3}.
\end{eqnarray}
The same argument applies to the interval $(1-C_0\delta_n,1)$.
Consequently, by (\ref{eq-oszacowania-1}), (\ref{eq-oszacowania}),
(\ref{eq-oszacowania-2}), (\ref{eq-oszacowania-3}) and comparing
$(y(1-y))^{1+\mu}$ with $(y(1-y))^{\mu}\psi_1(y)$, the statement of
Lemma \ref{lem.1} remains true when multiplying by $\psi_2(y)$. The
same holds true for Lemmas \ref{lem.3}, \ref{lem.2} and Corollary
\ref{cor.1}. Consequently,
Theorem \ref{th-reduction-quantiles-weighted} is proven.\\

The optimality of the bound in Theorem
\ref{th-reduction-quantiles-general} follows from Theorem
\ref{th-reduction-quantiles-weighted} in the same way we proved
optimality in Theorem \ref{cor-reduction-quantiles}. \koniec

\subsection{Proof of Theorem \ref{th-BK-general}}
Let $\beta<\frac{3}{4}$. Applying a third order Taylor expansion to
$f(Q(y))q_n(y)$, one has
\begin{eqnarray*}\lefteqn{
\hspace*{-1cm}|u_n(y)-f(Q(y))q_n(y)+\sigma_{n,1}n^{-1}\frac{f^{(1)}(Q(y))}{2f^2(Q(y))}u_n^2(y)|}\\
&=&\sigma_{n,1}^2n^{-2}\frac{f(Q(y))(y(1-y))^{3/2}}{6}Q^{(3)}(\theta)\sigma_{n,1}^{-3}n^{3}\frac{|y-U_n(y)|^3}{(y(1-y))^{3/2}}.
\end{eqnarray*}
We have
$$Q^{(3)}(y)=\frac{f^{(2)}(Q(y))}{f^4(Q(y))}-\frac{3(f^{(1)}(Q(y)))^2}{f^5(Q(y))}.$$
By the same argument as the one leading to (\ref{eq-oszacowania}),
it suffices to control the second term. We have
\begin{eqnarray*}\lefteqn{
(y(1-y))^{1/2}f(Q(y))\frac{(f^{(1)}(Q(\theta)))^2}{f^5(Q(\theta))}(y(1-y))^{3/2}}\\
&=&\frac{f(Q(y))}{f(Q(\theta))}\left(\frac{f^{(1)}(Q(\theta))}{f^2(Q(\theta))}\theta(1-\theta)\right)^2\left(\frac{y(1-y)}{\theta(1-\theta)}\right)^2.
\end{eqnarray*}
Under (CsR3(i)), (CsR3(ii)), in view of \cite[Lemma
1]{CsorgoRevesz1978} one has
\begin{equation}\label{CR-regularvariation}
\frac{f(Q(y))}{f(Q(\theta))}\le \left\{\frac{y\vee \theta}{y\wedge
\theta}\times \frac{1-(y\wedge \theta)}{1-(y\vee \theta
)}\right\}^{\gamma}.
\end{equation}
From this, (\ref{step.4}), (\ref{eq-bound-on-y-theta}) and Lemma
\ref{lem-weighted-LIL-quantile}, as $n\to\infty$, one concludes
\begin{eqnarray}\label{eq-modified-BK}
\lefteqn{\hspace*{-4cm}\sup_{y\in
(C_0\delta_n,1-C_0\delta_n)}(y(1-y))^{1/2}|u_n(y)-f(Q(y))q_n(y)+\frac{\sigma_{n,1}}{n}\frac{f(Q(y))f^{(1)}(Q(y))}{2f^3(Q(y))}u_n^2(y)|}\nonumber\\
&=&O_{\assubs}(\sigma_{n,1}^2n^{-2}(\log\log n)^{3/2}).
\end{eqnarray}
Next, taking Taylor expansion for $(\tilde V_{n,1}(y)-\tilde
V_{n,1}(U_n(y)))$, one obtains
\begin{eqnarray*}
\lefteqn{\sigma_{n,1}^{-1}(\tilde V_{n,1}(y)-\tilde
V_{n,1}(U_n(y)))=}\\
&= &
\sigma_{n,1}^{-1}\frac{f^{(1)}(Q(y))}{f(Q(y))}(y-U_n(y))\sum_{i=1}^nX_i+\sigma_{n,1}^{-1}(f\circ
Q)^{(2)}(\theta)(y-U_n(y))^2\sum_{i=1}^nX_i.
\end{eqnarray*}
Like in (\ref{eq-oszacowania}), as $n\to\infty$,
\begin{eqnarray}\label{step.300}
\lefteqn{\sup_{y\in
(C_0\delta_n,1-C_0\delta_n)}(y(1-y))^{\mu}\sigma_{n,1}^{-1}(f\circ
Q)^{(2)}(\theta)(y-U_n(y))^2\left|\sum_{i=1}^nX_i\right|}\\
&=&O_{\assubs}(\sigma_{n,1}^2n^{-2}(\log\log n)^{3/2}).\nonumber
\end{eqnarray}
If $\beta\ge \frac{3}{4}$, (\ref{eq-modified-BK}) and
(\ref{step.300}) remain valid if one replaces the weight functions
with $(y(1-y))^2$.

Thus,
\begin{eqnarray*}
\lefteqn{\sup_{y\in (C_0\delta_n,1-C_0\delta_n)}\psi_3(y)\left|\alpha_n(y)-f(Q(y))q_n(y)-\sigma_{n,1}^{-1}n^{-1}\frac{f^{(1)}(Q(y))}{2}\left(\sum_{i=1}^nX_i\right)^2\right|}\\
&\le &\mbox{\rm left hand side of } (\ref{eq-modified-BK})\\&&+
\sup_{y\in (C_0\delta_n,1-C_0\delta_n)}\psi_3(y)|\tilde
R_n(y)+\sigma_{n,1}^{-1}(\tilde
V_{n,1}(y)-\tilde V_{n,1}(U_n(y)))|\\
&&+\sup_{y\in
(C_0\delta_n,1-C_0\delta_n)}\psi_3(y)\left|\sigma_{n,1}^{-1}(\tilde
V_{n,1}(y)-\tilde
V_{n,1}(U_n(y)))\right.\\
&&\left.+\sigma_{n,1}n^{-1}\frac{f^{(1)}(Q(y))}{2f^2(Q(y))}u_n^2(y)+\sigma_{n,1}^{-1}n^{-1}\frac{f^{(1)}(Q(y))}{2}\left(\sum_{i=1}^nX_i\right)^2\right|\\
&= & O_{\assubs}(\sigma_{n,2}n^{-2}(\log\log
n)^{3/2})+O_{\assubs}(d_{n,p})+O_{\assubs}(\sigma_{n,1}^{-1}b_n(\log n)^{1/2})\\
&&+ \sigma_{n,1}n^{-1}\sup_{y\in
(C_0\delta_n,1-C_0\delta_n)}\psi_3(y)\left|\frac{f^{(1)}(Q(y))}{2}\right|\left|\frac{u_n(y)+\sigma_{n,1}^{-1}\sum_{i=1}^nX_i}{f(Q(y))}\right|^2\\
&&+O_{\assubs}(\sigma_{n,1}^2n^{-2}(\log\log n)^{1/2})
\end{eqnarray*}
by (\ref{eq-modified-BK}), (\ref{step.300}) and
(\ref{BahadurKiefer-approximation-step}) together with Lemmas
\ref{lem.2}, \ref{cor.1}. Moreover, by (CsR3) and via Theorem
\ref{th-reduction-quantiles-general} the bound is of the
order $O_{\assubs}(c_n)+O_{\assubs}(d_{n,p})$, as $n\to\infty$.\\

Further, as $n\to\infty$, $$\sup_{y\in
(0,C_0\delta_n)}\psi_3(y)\sigma_{n,1}^{-1}n^{-1}\left(\sum_{i=1}^nX_i\right)^2=O_{\assubs}(\delta_n^{1+\mu}\sigma_{n,1}n^{-1}(\log\log
n))=O_{\assubs}(c_n),$$ and $\sup_{y\in
(0,C_0\delta_n)}\psi_3(y)|\alpha_n(y)|=O_{\assubs}(c_n)$.

Next, having tail monotonicity assumption (CsR4) we may proceed as
in \cite{CsorgoRevesz1978}. Let $(k-1)/n<y\le k/n$. If $U_{k:n}\ge
y$, then $$\sup_{y\in (0, C_0\delta_n)}\psi_3(y)|f(Q(y))q_n(y)|\le
\sup_{y\in (0,
C_0\delta_n)}\psi_3(y)|u_n(y)|=O_{\assubs}(\delta_n^{(1+\mu)})=O_{\assubs}(c_n).$$
Further, if $U_{k:n}\le y$, then
$$\sup_{y\in (0,
C_0\delta_n)}\psi_3(y)|f(Q(y))q_n(y)|\le C
\sigma_{n,1}^{-1}n\sup_{y\in (0,
C_0\delta_n)}y(y(1-y))^{1+\mu}\log(\delta_n/U_{k:n})
$$
for $\gamma_1=1$. Now,
\begin{equation}\label{eq-orderstatistics}
P(U_{1:n}\le n^{-2}(\log n)^{-3/2})\le \sum_{i=1}^nP(U_i\le
n^{-2}(\log n)^{-3/2})\le n^{-1}(\log n)^{-3/2}.
\end{equation}
Consequently, via the Borel-Cantelli Lemma, as $n\to\infty$,
$U_{k:n}^{-1}=o_{\assubs}(n^2(\log n)^{3/2})$. Therefore,
\begin{equation}\label{eq.around.zero}
\sup_{y\in (0,
C_0\delta_n)}\psi_3(y)|f(Q(y))q_n(y)|=O_{\assubs}(c_n)
\end{equation}
follows for $\gamma_1=1$. \koniec
\subsection{Proof of Proposition
\ref{th-approx-unifquantile-generalquantile}} We follow lines of the
proof from \cite[Theorem 3]{CsorgoRevesz1978}. In view of Lemma
\ref{lem-weighted-LIL-quantile} and the Taylor expansion of
$f(Q(y))q_n(y)$, the approximation is valid on
$(C_0\delta_n,1-C_0\delta_n)$, provided $\beta<\frac{3}{4}$. For
$\beta\ge\frac{3}{4}$ it remains true by the choice of $\psi_4(y)$.

Having tail monotonicity assumption (CsR4), let $(k-1)/n<y\le k/n$.
If $U_{k:n}\ge y$, then (cf. (3.13) in \cite{CsorgoRevesz1978})
$$\sup_{y\in (0,
C_0\delta_n)}\psi_4(y)|f(Q(y))q_n(y)|\le \sup_{y\in (0,
C_0\delta_n)}\psi_4(y)|u_n(y)|=O_{\assubs}(a_n)$$ from
(\ref{eq-small-increments-quantile}) if $\beta<\frac{3}{4}$, and by
the choice of $\psi_4(y)$ if $\beta\ge \frac{3}{4}$.

If $U_{k:n}\le y$ and $\beta\in (\frac{1}{2},1)$, then for
$\gamma=1$, as $n\to\infty$, $$\sup_{y\in (0,
C_0\delta_n)}|f(Q(y))q_n(y)|=O_{\assubs}(\sigma_{n,1}n^{-1}\ell
(n))$$ by (\ref{eq-orderstatistics}). Moreover, as in
(\ref{eq-orderstatistics}), $U_{k:n}^{-1}=o_{P}(n(\log n)^{3/2})$,
as $n\to\infty$. Therefore, for $\gamma>1$, as $n\to\infty$,
$$\sup_{y\in (0,
C_0\delta_n)}|f(Q(y))q_n(y)|=O_{P}(\sigma_{n,1}n^{-1}\ell
(n)).$$\koniec

\end{document}